\documentclass[11pt]{article}
\usepackage{}

\usepackage[total={6in, 8.5in}]{geometry}
\usepackage{amsfonts}
\usepackage{amsthm}
\usepackage{amssymb}
\usepackage{amsmath}
\usepackage{enumerate}
\usepackage[
pdfauthor={ESYZ},
pdftitle={Toughness and spanning trees in K4mf graphs},
pdfstartview=XYZ,
bookmarks=true,
colorlinks=true,
linkcolor=blue,
urlcolor=blue,
citecolor=blue,
bookmarks=false,
linktocpage=true,
hyperindex=true
]{hyperref}

\makeatletter
\newcommand{\setword}[2]{%
	\phantomsection
	#1\def\@currentlabel{\unexpanded{#1}}\label{#2}%
}
\makeatother

\makeatletter
\renewcommand*\env@matrix[1][*\c@MaxMatrixCols c]{%
	\hskip -\arraycolsep
	\let\@ifnextchar\new@ifnextchar
	\array{#1}}
\makeatother

\usepackage{graphicx}
\usepackage[natural]{xcolor}
\usepackage{tikz}
\usepackage{tikz-3dplot}
\usetikzlibrary{shapes}
\usetikzlibrary{arrows,decorations.pathmorphing,backgrounds,positioning,fit,petri,automata}
\usetikzlibrary {positioning}
\usetikzlibrary{arrows}

\usepackage{pgf,tikz,pgfplots}

\usepackage{mathrsfs}

\usetikzlibrary{arrows}

\usepackage{xcolor}
\long\def\ignore#1{}

\let\oldi\ignore

\newcommand{\D}{\Delta}

 \oldi{
\newtheorem{THM}{\textbf{Theorem}}[section]
\newtheorem{THMs}{\textbf{Theorem}}[section]
\newtheorem{DEF}[THM]{\textbf{Definition}}[section]
\newtheorem{LEM}[THM]{\textbf{Lemma}}
\newtheorem{CON}[THM]{\textbf{Conjecture}}
\newtheorem{PROP}[THM]{\textbf{Proposition}}
\newtheorem{COR}[THM]{\textbf{Corollary}}
\newtheorem{CORs}{\textbf{Corollary}}[section]
\newtheorem{PRO}[THM]{\textbf{Problem}}
\newcommand{\pf}{\textbf{Proof}.\quad}

\newtheorem{FAC}{\textbf{Fact}}
\newtheorem{REM}{\textbf{Remark}}
\newtheorem{OPR}{\textbf{Operation}}
\newtheorem{CLA}{\textbf{Claim}}[section]
\setcounter{CLA}{1}
 }

\newtheorem{THM}{Theorem}[section]

\newtheorem{DEF}[THM]{Definition}
\newtheorem{LEM}[THM]{Lemma}
\newtheorem{CON}[THM]{Conjecture}

\newtheorem{CLA}{Claim}[section]
\newcommand{\pf}{\textbf{Proof}.\quad}

\newtheorem*{THM2}{\textbf{Theorem 1.4}}
\newtheorem*{THM3}{\textbf{Theorem 1.5}}
\newtheorem*{THM4}{\textbf{Theorem 2.5}}

\usepackage{thmtools}
\usepackage{thm-restate}

\newcommand{\CC}{\mathcal{C}}

\newcommand{\phiv}{\varphi}

\newcommand{\pbar}{\overline{\varphi}}

\begin{document}
\title{Overfullness of  critical class 2 graphs with a small core degree}

\author{%
 Yan Cao
 \qquad Guantao Chen\thanks{This work was supported in part by NSF grant DMS-1855716.}\\
  Department of Mathematics and Statistics, \\
  Georgia State University, Atlanta, GA 30302, USA\\
   \texttt{ycao17@gsu.edu}
  \qquad
   \texttt{gchen@gsu.edu}%
 \and 
 Songling Shan\\
 Department of Mathematics, \\
 Illinois State Univeristy, Normal, IL 61790, USA \\
 \texttt{sshan12@ilstu.edu}
} 

\date{\today}
\maketitle

 \begin{abstract}
 Let $G$ be a simple graph, and let  $n$, $\Delta(G)$ and $\chi' (G)$ be the order, the maximum degree and the chromatic index of $G$, respectively.   We call $G$ 
  \emph{overfull} if $|E(G)|/\lfloor n/2\rfloor > \D(G)$, and {\it critical} if $\chi'(H) < \chi'(G)$ for every proper subgraph $H$ of $G$.  Clearly, if $G$ is overfull then $\chi'(G) = \Delta(G)+1$. The \emph{core} of $G$, denoted by $G_{\D}$, is the subgraph of $G$ induced by all its maximum degree vertices.  Hilton and Zhao conjectured that  for any critical class 2 graph $G$ with $\D(G) \ge 4$, if the maximum degree of $G_{\D}$ is at most two, then $G$ is overfull, which in turn gives $\D(G) > n/2 +1$.  We show that for any critical class 2 graph $G$, if the  minimum degree of $G_{\D}$ is at most two  and $\D(G) > n/2 +1$, then $G$ is overfull.   
 

 \end{abstract}


\section{Introduction}
We will mainly adopt the notation from the book~\cite{StiebSTF-Book}. Graphs in this paper are simple, i.e., finite, undirected, without loops or  multiple edges. Let $G$ be a graph.  A {\it $k$-edge-coloring\/} of $G$ is a map $\varphi$:  $E(G) \rightarrow \{1, 2, \dots, k\}$   that assigns to every edge $e$ of $G$ a color $\varphi(e) \in \{1, 2, \dots, k\}$  such that  no two adjacent edges receive the same color.  Denote by $\CC^k(G)$ the set of all $k$-edge-colorings of $G$. 
The {\it chromatic index\/} $\chi'(G)$ is  the least integer $k\ge 0$ such that $\CC^k(G) \ne \emptyset$.  Denote by 
$\D(G)$ the maximum degree of $G$. 
In 1960's, Vizing~\cite{Vizing-2-classes} and, independently,  Gupta~\cite{Gupta-67} proved that $\D(G) \le \chi'(G) \le \D(G)+1$. 
This leads to a natural classification of graphs. Following Fiorini and Wilson~\cite{fw},  we say a graph $G$ is of {\it class 1} if $\chi'(G) = \D(G)$ and of \emph{class 2} if $\chi'(G) = \D(G)+1$.
Holyer~\cite{Holyer} showed that it is NP-complete to determine whether an arbitrary graph is of class 1.

A graph G is \emph{critical} if $\chi'(H)<\chi'(G)$ for every proper subgraph $H$ of $G$.
In investigating the classification problem, critical graphs are of particular interest. 
Critical graphs of class 2 have rather more structure than arbitrary graphs of class two, and it follows from Vizing's theorem that every graph of
class 2 contains a critical graph of class 2 with the same maximum degree as a subgraph. In  this paper, we call a critical class 2 graph  \emph{$\D$-critical} if
 $\D(G) = \D$.

 Since every matching of  $G$ has at most 
 $\lfloor |V(G)|/2\rfloor$ edges,  $\chi'(G) \ge |E(G)|/ \lfloor |V(G)|/2\rfloor$. 
 A graph $G$ is {\em overfull} if $|E(G)| / \lfloor |V(G)|/2\rfloor > \D(G)$.  Clearly,  if $G$ is overfull then $\chi'(G)  = \D(G) +1$, and so $G$ is of class 2. Applying Edmonds' matching polytope theorem, Seymour~\cite{seymour79}  showed  that whether a graph containing an overfull subgraph can be determined in polynomial time.  A number of outstanding conjectures listed in {\it Twenty Pretty Edge Coloring Conjectures} in~\cite{StiebSTF-Book} lie in deciding when a $\D$-critical graph is overfull.   
  
 The \emph{core} of a graph $G$, denoted by $G_{\D}$, is the subgraph induced by all its maximum degree vertices.
 Vizing~\cite{Vizing-2-classes} proved that if $G_\Delta$ has at most two vertices then $G$ 
is class 1. Fournier~\cite{MR0349458} generalized Vizing's result by showing that 
if $G_\Delta$ is acyclic then $G$ is class 1. 
Thus a necessary condition for a graph to be class 2 is to have a core
that contains cycles. A long-standing conjecture of Hilton and Zhao~\cite{MR1395947} claims that for a $\D$-critical graph $G$ with $\D\ge 4$, 
if the maximum degree of $G_{\D}$ is at most two, then $G$ is overfull. We~\cite{HZ}, along with Guangming Jing, recently confirmed this conjecture, which in turn implies  $\D(G) > n/2 +1$, where $n =|V(G)|$ is the order of $G$. In this paper, by imposing a condition on the maximum degree of $G$, we relax the condition $\Delta(G_{\D}) \le 2$, and show a result analogous to 
the Hilton-Zhao Conjecture as below. 

\begin{THM}\label{thm:main}
Let $G$ be a $\Delta$-critical graph of order $n$.  If $\delta(G_\D) \le 2$
and $\Delta(G) > n/2+1$,  
then $G$ is overfull. 
\end{THM}

Implicitly, Theorem~\ref{thm:main} is much stronger than the Hilton-Zhao Conjecture, but we don't have a direct proof for that. 
A graph $G$ is said to be {\it just overfull} if $|E(G)| = \D(G) \lfloor \frac 12 |V(G)|\rfloor +1$. 
 We  hope that the new edge coloring techniques we introduced in our proof may shed some light on attacking the Just Overfull Conjecture -- Conjecture 4.23 (page 72) in ~\cite{StiebSTF-Book}. 
 
 \begin{CON}\label{Conj-just}
Let $G$ be a $\D$-critical graph of order $n$.  If $\D(G) \ge n/2$, then $G$ is just overfull.  
\end{CON}

Chetwynd and  Hilton in 1986~\cite{MR848854,MR975994} made a much stronger conjecture,  commonly referred to as  {\it the Overfull Conjecture} that for a $\D$-critical graph of order $n$, if $\D(G) > n/3$ then $G$ is overfull. Except some very special results~\cite{MR975994,comm,MR2082738}, the Overfull Conjecture seems untouchable with current edge-coloring techniques.  

Let $G$ be a graph and $H \subseteq  G$ be a subgraph.  For $v\in V(G)$, $N(v)$  is the set of neighbors of $v$
in $G$ and $d(v)=|N(v)|$ is the degree of $v$ in $G$. Let 
 $N_H(v)=N(v)\cap V(H)$ and $d_H(v)=|N_H(v)|$.  More generally, for a subset $S\subseteq V(G)$, let $N_H(S) = \cup_{v\in S} N_H(v)$ be the neighborhood of  $S$ in $G$ that is contained in  $V(H)$.   
  For a nonnegative integer $k$, 
 a {\it $k$-vertex} is a vertex of degree $k$. 
We denote by $V_k$ and $N_k(v)$  the set of all  $k$-vertices,  repetitively, in  $V(G)$ and  $N(v)$. 
Let $N[v]=N(v)\cup \{v\}$ and $N_k[v]=N_k(v)\cup \{v\}$. 
For convenience, for any nonnegative integers $i$ and $j$, let $[i, j] = \{i, i+1, \dots, j\}$. 

A vertex $v$ is called {\it light} if it is adjacent to at most two $\D(G)$-vertices, i.e., $d_{G_\D}(v) \le 2$. 
An edge $e$ of a graph $G$ is {\it critical} if 
$\chi'(G -e) < \chi'(G)$. Clearly, if $G$ is $\D$-critical then every edge of $G$ is critical. In a $\D$-critical graph, for a light vertex,  what can we say  about its neighbors? 
The following lemma reveals  some of their properties. 

\begin{LEM}[Vizing's Adjacency Lemma (VAL)] Let $G$ be a class 2 graph with maximum degree $\Delta$. If $e=xy$ is a critical edge of $G$, then $x$ is adjacent to at least $\Delta-d(y)+1$ $\Delta$-vertices from $V(G)\setminus \{y\}$.
	\label{thm:val}
\end{LEM}

Let $G$ be a $\Delta$-critical graph and $r$ be a  light vertex of $G$.  We claim  $d(s) \ge \Delta-1$ for every $s\in N(r)$.   
Otherwise, by VAL, $r$ is adjacent to at least $\Delta-d(s)+1\ge 3$ vertices of degree $\D$,
giving a contradiction. Consequently, we have $d(s) = \D -1$ or $d(s) = \D$.  We also see that $r$ must be adjacent to exactly 
two $\D$-vertices if $r$ is light. These facts will be frequently used throughout this paper. 

Theorem~\ref{thm:main} is a consequence of the following three technical results.  

\begin{THM}\label{thm:s1-adj}
	Let $G$ be a  class 2 graph with maximum degree $\Delta$, and 
	$r\in V(G)$  be a light $\D$-vertex. For any  $s\in N(r)$ with $d(s)<\Delta$, if  $rs$ is a critical edge 
	of $G$,  then all vertices in $N(s)\backslash N(r)$ are $\D$-vertices. 
\end{THM}

\begin{THM}\label{thm:longk}
	Let $G$ be a  class 2 graph with maximum degree $\Delta$,  $r\in V(G)$ be a light $\D$-vertex  and $s\in N_{\D-1}(r)$ such that $rs$ is a critical edge.  For every $x\in V(G)\setminus N[r]$, if  $d(x) \le \D -3$,  then $N(x)\cap N(s) \subseteq N(r)\setminus N_\Delta(r)$. 
\end{THM}

\begin{THM}\label{thm:longk2}
Let $G$ be a  $\Delta$-critical graph of order $n$.  If  $\Delta >n /2+1$ and $\delta(G_{\D}) \le 2$,  
then $n$ is odd. 
\end{THM}

\vskip .2in
\proof [{\bf Proof of Theorem~\ref{thm:main}}]
Let $G$ be a $\D$-critical graph of order $n$ such that $\delta(G_{\D}) \le 2$ and $\D > n/2+1$. 
By Theorem~\ref{thm:longk2}, $n$ is odd.
Let $r$ be a light $\D$-vertex of $G$ and $S\subseteq N_G(r)$ such that $d(s) < \D$ for each $s\in S$.  By the remark immediately  after Lemma~\ref{thm:val},  we have $d(s) = \D-1$ for every $s\in S$. Since $|S|=\D-2$, we have $2|E(G)|\le n \D-(\D-2)$. Thus to show $2|E(G)| \ge (n-1)\D+2$ (i.e., $G$ is overfull), 
we only need to  show that 
all vertices in $V(G)\setminus S$ are $\D$-vertices.  

Assume to the contrary that there exists $x\in V(G)\setminus S$ such that  $d(x) \le \Delta-1$.  Since every vertex in 
$N[r]\setminus S$ is a $\D$-vertex,  we have  $x\notin N[r]$.  Since $n$ odd, 
$\D > n/2 +1$ implies  $\D \ge (n+1)/2+1$. 
We first suppose that  $d(x)\ge \Delta-2$,  i.e.,  $|N(x)|  \ge (n-1)/2$. 
Since  $|S|=\Delta-2\ge  (n-1)/2$ and $r\not\in N(x)$, 
we conclude that $N(x)\cap S \ne \emptyset$. Let $s\in N(x)\cap S$. Since $G$ is $\Delta$-critical, $rs$ is a critical edge of $G$. 
Applying Theorem~\ref{thm:s1-adj},  we get $d(x) =\D$, a contradiction.   Thus $d(x)\le \Delta-3$. 
Since $G$ is $\D$-critical,  $x$ has a neighbor $u$ with degree $\D$.
 As $\D \ge (n+1)/2+1$ and $|S| = \D -2$, we find a vertex $s\in N(u)\cap S$. Thus $u\in N(x)\cap N(s)$.  
 Since $d(u)=\Delta$ and $d(x)\le \Delta-3$, $u\not\in N(r)\setminus N_\Delta(r)$.
 Again, $rs$ is a critical edge of $G$ as $G$ is $\D$-critical.
Applying the contrapositive statement of  Theorem~\ref{thm:longk}, we get $d(x) \ge \D -2$, which gives a contradiction. 
\qed

Theorems~\ref{thm:s1-adj} to \ref{thm:longk2} study some structural properties
 of vertices outside the neighborhood of a light vertex. 
The study of structural properties of vertices beyond a given neighborhood plays a key role in our proof, and we believe that the technique may be useful on tackling other edge-coloring problems involving overfull properties.

\section{Preliminaries}\label{lemma}
This section is divided into three subsections. We first give some basic notation and terminologies, then define a slightly modified and specific Vizing fan centering at a light vertex, and finally we investigate some properties of 
a $\Delta$-coloring  around a light vertex. 

\subsection{Basic notation and terminologies}

Let $G$ be a graph with maximum degree $\D$, and  let $e\in E(G)$ and $\varphi \in \CC^{\D}(G-e)$.  When we apply some definitions later, we may drop the phrase ``w.r.t. $\varphi$"  or surpass the coloring symbol $\varphi$ whenever the edge-coloring $\varphi$ is clearly understood.

For a vertex  $v\in V(G)$, define the two color sets
\[
\varphi(v)=\{\varphi(f)\,:\, \text{$f\ne e$ is incident to $v$}\} \quad \mbox{ and}\quad \pbar(v)=[1, \D] \setminus\varphi(v).
\]
We call $\varphi(v)$ the set of colors \emph{present} at $v$ and $\pbar(v)$
the set of colors \emph{missing} at $v$. 
If $|\pbar(v)|=1$, we will also use $\pbar(v)$ to denote the  color  missing at $v$.

For a vertex set $X\subseteq V(G)$,  define  $\pbar(X)=\bigcup _{v\in X} \pbar(v)$ to be the set of missing colors of $X$. 
The set $X$ is  \emph{elementary} w.r.t. $\varphi$  or simply \emph{$\varphi$-elementary} if $\pbar(u)\cap \pbar(v)=\emptyset$
for every two distinct vertices $u,v\in X$.   

For a color $\alpha$, the edge set $E_{\alpha} = \{ e\in E(G)\, |\, \varphi(e) = \alpha\}$ is called a {\it color class}. Clearly, 
$E_{\alpha}$ is a {\it matching} of $G$ (possibly empty). 
For two distinct colors $\alpha,\beta$,  the  subgraph of $G$
induced by $E_{\alpha}\cup E_{\beta}$ is a union of disjoint 
paths and  even cycles, which are referred to as   \emph{$(\alpha,\beta)$-chains} of $G$
w.r.t. $\varphi$.  For a vertex $v$, let $C_v(\alpha, \beta, \varphi)$ denote the unique $(\alpha, \beta)$-chain 
containing $v$.  
If $C_v(\alpha, \beta, \varphi)$ is a path, we just write it as $P_v(\alpha, \beta, \varphi)$. The latter is commonly used when we know that $|\pbar(v)\cap \{\alpha,\beta\}|=1$.  If we interchange the colors $\alpha$ and $\beta$
on an $(\alpha,\beta)$-chain $C$ of $G$, we briefly say that the new coloring is obtained from $\varphi$ by an 
{\it $(\alpha,\beta)$-swap} on $C$, and we write it as  $\varphi/C$. 
This operation is called a \emph{Kempe change}.  If $C=uv$
is just an edge, the notation  $\mathit{uv: \alpha\rightarrow \beta}$  means to recolor  the edge  $uv$ that is colored by $\alpha$ using the color $\beta$.


%


 Suppose that $\alpha, \beta, \gamma$ are three colors such that  $\alpha\in \pbar(x)$ and  $\beta,\gamma\in \varphi(x)$. An $\mathit{(\alpha,\beta)-(\beta,\gamma)}$
\emph{swap at $x$}  consists of two operations:  first swaps colors on $P_x(\alpha,\beta, \varphi)$ to get a new coloring $\varphi'$, and then swaps
colors on $P_x(\beta,\gamma, \varphi')$. 
 When $\beta =\alpha$,  an $(\alpha,\alpha)$-swap is just a vacuous recoloring operation. 
 
 For a given path $P$, a vertex $u$ and an edge $uv$, we write $u\in P$ and $uv\in P$ for $u\in V(P)$ and $uv\in E(P)$, respectively. 
Suppose $x\in P$. For two vertices $u,v\in P_x(\alpha,\beta, \varphi)$, if  $u$ lies between $x$ and $v$, 
then we say that $P_x(\alpha,\beta, \varphi)$ \emph{meets $u$ before $v$}.

\subsection{Modified Vizing fans and Kierstead paths }

The fan argument was introduced by Vizing~\cite{Vizing64,vizing-2factor} in his classic results on the upper bounds of chromatic indices.  We will use multifan, a generalized version of Vizing fan, given by Stiebitz et al.~\cite{StiebSTF-Book}, in our proof. To simplify the arguments, we will not include maximum degree vertices in our fans except the center vertex. 

\begin{DEF} Let  $G$ be a graph with maximum degree $\D$.    For an edge $e=rs_1\in E(G)$ and  a coloring $\varphi\in \CC^{\D}(G-e)$, 
a \emph{multifan} centered at $r$ w.r.t. $e$ and $\varphi$
is a sequence $F_\varphi(r,s_1:s_p)=(r, rs_1, s_1, rs_2, s_2, \ldots, rs_p, s_p)$ with $p\geq 1$ consisting of  distinct vertices $r, s_1,s_2, \ldots , s_p$ and edges $rs_1, rs_2,\ldots, rs_p$ satisfying   the following condition:
\begin{enumerate}  [{\em (F1)} ]
	\item For every edge $rs_i$ with $i\in [2, p]$,  there exists  $j\in [1, i-1]$ such that 
	$\varphi(rs_i)\in \pbar(s_j)$, 
\end{enumerate}
and none of $s_1, \dots, s_p$ is a $\D$-vertex. 
\end{DEF}
We will simply denote a multifan  $F_\varphi(r,s_1: s_{p})$ by $F$ if we do not need to emphasize the center $r$,  and the non-center
starting and  ending vertices. We also notice that if $F_\varphi(r,s_1: s_{p})$ is a multifan, then for any integer
$p^*\in [1,  p]$, $F_\varphi(r,s_1: s_{p^*})$
is also a multifan.  The following result regarding a multifan can be found in \cite[Theorem~2.1]{StiebSTF-Book}.

\begin{LEM}
	\label{thm:vizing-fan1}
	Let $G$ be a class 2 graph, $e=rs_1$ be a critical edge and $\varphi\in \CC^\Delta(G-e)$. 
	If $F_\varphi(r,s_1:s_p)$  is a multifan w.r.t. $e$ and $\varphi$,   then  $V(F)$ is $\varphi$-elementary. 
\end{LEM}

Suppose that $e=rs_1$ is a critical edge of a class 2 graph $G$ and  $F_\varphi(r,s_1:s_p)$  is a multifan w.r.t. $e$ and a coloring $\varphi\in \CC^{\D}(G-e)$.  Given a color $\alpha\in \pbar(s_1)$, we call a vertex $s_{\ell}$ with $\ell\in [2,p]$ an \emph{$\alpha$-inducing vertex} if 
there exists a subsequence $(s_{\ell_1},s_{\ell_2}, \ldots, s_{\ell_k})$ terminated at $s_{\ell_k} = s_\ell$ such that $
\varphi(rs_{\ell_1})= \alpha\in \pbar(s_1)$ and for each $i\in [2,k]$,
 $\varphi(rs_{\ell_i})\in \pbar(s_{\ell_{i-1}})$.   We also call the above sequence an $\alpha$-sequence, and a color
 $\beta\in \pbar(s_{\ell})$ an $\alpha$-inducing color or a color induced by $\alpha$.   For convention, $\alpha$ itself is also called an $\alpha$-inducing color. Since $V(F)$ is elementary, every color in $\pbar(F\backslash \{r\})$ is induced by a color in $\pbar(s_1)$. 
 
 As a consequence of Lemma~\ref{thm:vizing-fan1}, we have the following linkage properties of vertices in a multifan. The original proof can  be found in~\cite{HZ}. 

\begin{LEM}
	\label{thm:vizing-fan2}
	Let $G$ be a class 2 graph, $e=rs_1$ be a critical edge and  $\varphi\in \CC^\Delta(G-e)$. 
	Then, for every multifan  $F_\varphi(r,s_1:s_p)$,  the following three statements hold. 
	
		\begin{enumerate}[(a)]
	\item For any color $\gamma\in \pbar(r)$ and any color $\delta\in \pbar(s_i)$ with $i\in [1, p]$,  vertices $r$ and $s_i$ are $(\gamma,\delta)$-linked w.r.t. $\varphi$. \label{thm:vizing-fan1b}	
	
	\item For $i,j\in [1,p]$, if  two colors $\delta\in \pbar(s_i)$ and $\lambda\in \pbar(s_j)$ are induced by two different colors in $\pbar(s_1)$,  then the corresponding vertices $s_i$ and $s_j$ are $(\delta, \lambda)$-linked.  \label{thm:vizing-fan2-a}
		
	\item For $i,j\in [1,p]$,  suppose  two colors $\delta\in \pbar(s_i)$ and $\lambda\in \pbar(s_j)$ are induced by a same  color in $\pbar(s_1)$. 	  If $s_i$ and $s_j$ are not 
	$(\delta, \lambda)$-linked and $j>i$, then    
		$r\in P_{s_j}(\delta, \lambda, \varphi)$.  	\label{thm:vizing-fan2-b}
	\end{enumerate}
	
\end{LEM}
%

 Let $G$ be a class 2 graph, $r\in V(G)$ be a light vertex, $rs_1 \in E(G)$ be a critical edge and $\phiv\in \CC^{\D}(G-rs_1)$.  Let
 $F_\varphi(r,s_1:s_p)$ be a multifan w.r.t. $rs_1$ and $\phiv$.  By VAL, except  two $\D$-vertices, all other neighbors of $r$ are $(\D-1)$-vertices.  In particular, $d(s_i) =\D -1$ for all $i\in [1, p]$. Hence,  $|\pbar(s_1)| =2$ and $|\pbar(s_i)| =1$ for each $i\in [2, p]$.  Assume without loss of generality $\pbar(s_1) =\{2, \D\}$. Then,  all $2$-inducing vertices form a $2$-sequence and all $\D$-inducing vertices form a $\D$-sequence. By relabeling if necessary, we assume $(s_2, \dots, s_{\alpha})$ is a $2$-sequence  and $(s_{\alpha+1}, \dots, s_p)$ is a $\D$-sequence for some $\alpha \in [1,p]$, where we define $(s_2, \dots, s_{\alpha})$ to be the empty sequence if 
 $\alpha<2$. 
  We call a multifan {\it typical} at a light vertex $r$, denoted  by $F_\varphi(r, s_1:s_\alpha:s_\beta)$,  if $1\in \pbar(r)$, $\pbar(s_1)=\{2,\Delta\}$ and either $|V(F)|=2$ or $|V(F)|\ge 3$ 
 with the following two conditions. 
\begin{enumerate}[(1)]
	\item  $(s_2, \ldots, s_\alpha)$  is a $2$-inducing sequence 
	and $(s_{\alpha+1}, \ldots, s_{\beta})$  is a $\Delta$-inducing sequence of $F$. 
	\item For each $i\in [2,\beta]$, $\varphi(rs_i)=i$ and $\pbar(s_i)=i+1$ except when $i=\alpha+1\in [3,\beta]$.  In this case,	$\varphi(rs_{\alpha+1})=\Delta$ and $\pbar(s_{\alpha+1})=\alpha+2$.
\end{enumerate}
 
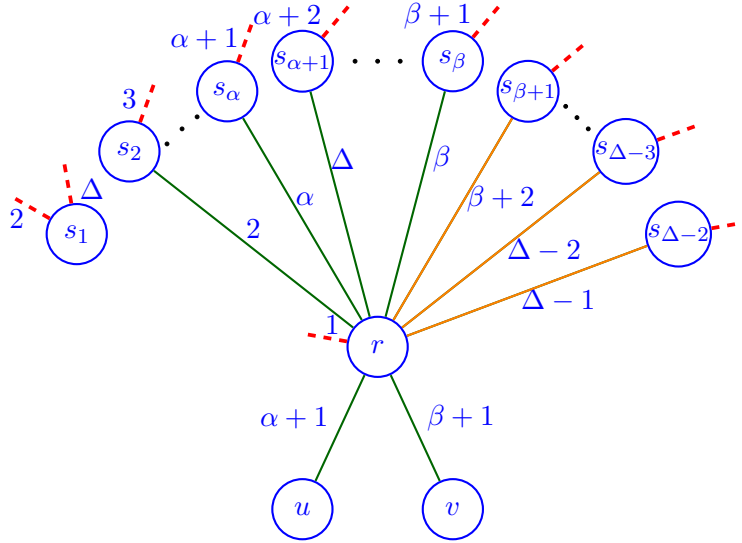
\begin{figure}[!htb]
	\begin{center}
		\begin{tikzpicture}[scale=1]
		
		{\tikzstyle{every node}=[draw ,circle,fill=white, minimum size=0.8cm,
			inner sep=0pt]
			\draw[blue,thick](0,-3) node (r)  {$r$};
			\draw [blue,thick](-4, -1.5) node (s1)  {$s_1$};
			\draw [blue,thick](-3.3, -0.4) node (s2)  {$s_2$};
			\draw[blue,thick] (-2, 0.4) node (sa)  {$s_\alpha$};
			\draw [blue,thick](-1, 0.8) node (sa2)  {$s_{\alpha+1}$};
			\draw [blue,thick](1, 0.8) node (sb)  {$s_{\beta}$};
			\draw [blue,thick](2, 0.4) node (sb2)  {$s_{\beta+1}$};
			\draw[blue,thick] (3.3, -0.4) node (sd1)  {$s_{\Delta-3}$};
			\draw [blue,thick](4, -1.5) node (sd2)  {$s_{\Delta-2}$};
			\draw [blue,thick](-1, -5.16) node (u)  {$u$};
			\draw [blue,thick](1, -5.16) node (v)  {$v$};
		}
		\path[draw,thick,black!60!green]
		(r) edge node[name=la,above,pos=0.5] {\color{blue}$2$} (s2)
		(r) edge node[name=la,pos=0.6] {\color{blue}\quad$\alpha$} (sa)
		(r) edge node[name=la,pos=0.7] {\color{blue}\quad$\Delta$} (sa2)
		(r) edge node[name=la,pos=0.7] {\color{blue}\quad$\beta$} (sb)
		(r) edge node[name=la,pos=0.6] {\color{blue}\qquad\,\,\,$\beta+2$} (sb2)
		(r) edge node[name=la,pos=0.5] {\color{blue}\quad \qquad$\Delta-2$} (sd1)
		(r) edge node[name=la,pos=0.4] {\color{blue}\qquad\quad\,\,\,\,\,$\Delta-1$} (sd2)
		(r) edge node[name=la,pos=0.4] {\color{blue}$\alpha+1$\qquad\quad\,\,\,} (u)
		(r) edge node[name=la,pos=0.4] {\color{blue}\qquad\quad\,\,\,$\beta+1$} (v);

		\draw[orange, thick] (r) --(sb2); 
		\draw [orange, thick](r) --(sd1); 
		\draw [orange, thick](r) --(sd2); 
		
		\draw[dashed, red, line width=0.5mm] (r)--++(170:1cm); 
		\draw[dashed, red, line width=0.5mm] (s1)--++(150:1cm); 
		\draw[dashed, red, line width=0.5mm] (s1)--++(100:1cm); 
		\draw[dashed, red, line width=0.5mm] (s2)--++(70:1cm); 
		\draw[dashed, red, line width=0.5mm] (sa)--++(70:1cm); 
		\draw[dashed, red, line width=0.5mm] (sa2)--++(50:1cm); 
		\draw[dashed, red, line width=0.5mm] (sb)--++(50:1cm); 
		\draw[dashed, red, line width=0.5mm] (sb2)--++(40:1cm); 
		\draw[dashed, red, line width=0.5mm] (sd1)--++(20:1cm); 
		\draw[dashed, red, line width=0.5mm] (sd2)--++(10:1cm); 
		
		\draw[blue] (-0.6, -2.7) node {$1$};  
		\draw[blue] (-4.8, -1.3) node {$2$};  
		\draw[blue] (-3.8, -0.9) node {$\Delta$};  
		\draw[blue] (-3.3, 0.3) node {$3$};  
		\draw[blue] (-2.3, 1.1) node {$\alpha+1$};  
		\draw[blue] (-1.2, 1.4) node {$\alpha+2$};  
		\draw[blue] (0.8, 1.4) node {$\beta+1$}; 
		
		{\tikzstyle{every node}=[draw ,circle,fill=black, minimum size=0.05cm,
			inner sep=0pt]
			\draw(-2.4,0.06) node (f1)  {};
			\draw(-2.6,-0.1) node (f1)  {};
			\draw(-2.8,-0.3) node (f1)  {};
			\draw(-0.3,0.8) node (f1)  {};
			\draw(0,0.8) node (f1)  {};
			\draw(0.3,0.8) node (f1)  {};
			\draw(2.5,0.2) node (f1)  {};
			\draw(2.65,0.05) node (f1)  {};
			\draw(2.8
			,-0.1) node (f1)  {};
		} 
		\end{tikzpicture}
	\end{center}
	\caption{A typical multifan $F_\varphi(r, s_1:s_\alpha:s_\beta)$ at a light vertex $r$, where a dashed line at a vertex indicates a color missing at the vertex.}
	\label{f11}
\end{figure}

A {\it typical multifan} at a light vertex $r$  is depicted in 
 Figure~\ref{f11}.  
 
 By relabeling if necessary, every multifan centered at a light vertex  $r$ 
 is corresponding to a typical multifan at $r$ on the same vertex set. Thus in this paper,
 we assume all mutilfans at $r$ are typical.

We close this subsection with Kierstead paths, which was introduced by Kierstead~\cite{K-path-84} in his work on edge-colorings of multigraphs. 
\begin{DEF} 
Let $G$ be a graph, $e=v_0v_1\in E(G)$, and  $\varphi\in \CC^{\D}(G-e)$.
A \emph{Kierstead path}  w.r.t. $e$ and $\varphi$
is a sequence $K=(v_0, v_0v_1, v_1, v_1v_2, v_2, \ldots, v_{p-1}, v_{p-1}v_p,  v_p)$ with $p\geq 1$ consisting of  distinct vertices $v_0,v_1, \ldots , v_p$ and  edges $v_0v_1, v_1v_2,\ldots, v_{p-1}v_p$ satisfying the following condition:
\begin{enumerate}[{\em (K1)}]
	\item For every edge $v_{i-1}v_i$ with $i\in [2,p]$,  there exists $j\in [1,i-2]$ such that 
	$\varphi(v_{i-1}v_i)\in \pbar(v_j)$. 
\end{enumerate}
\end{DEF}

Clearly a Kierstead path with at most three vertices is a multifan. So we consider 
Kierstead paths with four  vertices and restrict on its simple graph version. 
The following lemma was proved in Theorem 3.3 from~\cite{StiebSTF-Book}.  

\begin{LEM}[]\label{Lemma:kierstead path1}
	Let $G$ be a class 2 graph,
	$e=v_0v_1\in E(G)$ be a critical edge, and $K=(v_0, v_0v_1, v_1, v_1v_2,  v_2, v_2v_3, v_3)$ be a Kierstead path w.r.t. $e$ and  a coloring $\varphi\in \CC^{\D}(G-e)$.   If $\min\{d_G(v_2), d_G(v_3)\}<\Delta$,  then $V(K)$ is $\varphi$-elementary. 		
\end{LEM}

Let $G$ be a class 2 graph,  $e$ be a critical edge and $\phiv\in \CC(G-e)$.   Let $T$ be  a sequence of vertices  and edges of  $G$. We denote by \emph{$V(T)$}  
and \emph{$E(T)$} 
the set of vertices and the set of edges that are contained in $T$, respectively.  If $V(T)$ is $\varphi$-elementary,
then for a color  $\tau\in \pbar(V(T))$,  we denote by  $\mathit{\pbar^{-1}_T(\tau)}$ the  unique vertex  in $V(T)$ at which $\tau$ 
is missing.  
A coloring $\varphi' \in \CC^\Delta(G-e)$ is called {\it $T$-stable} w.r.t. $\varphi$ if 
$\pbar'(x)=\pbar(x)$
for every  vertex $x\in V(T)$  and 
$\varphi'(f)=\varphi(f)$
for every edge  $f\in E(T)$.   Clearly, $\varphi$ is 
$T$-stable w.r.t. itself. 
For simplicity, we write $\pbar(T)$
for $\pbar(V(T))$. 

Let $F=F_\varphi(r, s_1:s_{\alpha}:s_{\beta})$ be a typical multifan w.r.t. $e=rs_1$ and $\phiv\in \CC^\Delta(G-rs_1)$. By the  definition above, if $\phiv'$ is $F$-stable, then $F$ is also a typical multifan w.r.t. $e$ and $\phiv'$.  Let $\gamma,\delta\in [1,\D]$ be two colors and 
 $P$ be a $(\gamma, \delta)$-path. If $E(P)\cap E(F) = \emptyset$ and neither end-vertices of $P$ are in $V(F)$, then Kempe change $\phiv/P$ gives an $F$-stable coloring. Applying Lemma~\ref{thm:vizing-fan2}, we have the following results on stable coloring, which will be used heavily in our proofs. 

\begin{LEM}\label{Lem:Stable}
Let $G$ be a class 2 graph and $F=F_\varphi(r, s_1:s_{\alpha}:s_{\beta})$ be a typical multifan w.r.t. a light vertex $r$, critical edge $rs_1$, and a coloring $\phiv\in \CC^{\D}(G-rs_1)$. For any color $\gamma\in \pbar(F)$ and $x\notin V(F)$,  the following statements hold.
\begin{itemize}
\item the Kempe change $\phiv/P_x(1,\gamma, \phiv)$ gives an $F$-stable coloring provided $\pbar(x)\cap \{1, \gamma\} \ne \emptyset$. 

\item if $\gamma$ is $2$-inducing, then the Kempe change $\phiv/P_x(\gamma, \D, \phiv)$ gives an $F$-stable coloring provided 
$r\notin P_x(\gamma, \D, \phiv)$ and $\pbar(x)\cap \{\gamma, \D\} \ne \emptyset$;  and 

\item if $\gamma$ is $\D$-inducing,  then the Kempe change $\phiv/P_x(2,\gamma,  \phiv)$ gives an $F$-stable coloring provided 
$r\notin P_x(\gamma, 2, \phiv)$ and $\pbar(x)\cap \{\gamma, 2\} \ne \emptyset$. 
\end{itemize}

\end{LEM}

\subsection{$\tau$-sequence,  rotation, and shifting}
 Throughout this subsection, we assume that  $G$ is a class 2 graph, $r\in V(G)$ is a light vertex,    $e=rs_1\in E(G)$ is a critical edge of $G$ and $\varphi\in \CC^{\D}(G-e)$. We also assume  that  $N_{\D}(r) =\{u_1, u_2\}$ and $N_{\D -1}(r) =\{s_1, \dots, s_{q}\}$,  where $q=d(r)-2$.  Furthermore, we assume that $F=F_{\varphi}(r,s_1:s_{\alpha}:s_{\beta})$ is a typical multifan at $r$.

We call $F$ a {\it maximum mutlifan at $r$} if $|V(F)|$ is maximum over all colorings in $\CC^{\D}(G-e)$ and all multifans centered at $r$.
Clearly, if  $F$ is maximum, then  colors $\alpha+1$ and  $\beta+1$ are assigned to  edges $ru_1$ and $ru_2$, respectively,  i.e., 
$\alpha+1, \beta+1\notin \{\varphi(rs_{\beta+1}), \dots, \varphi(rs_q)\}$ (see Figure~\ref{f11}).

\begin{DEF} \label{Def-tauSeq} For a color  $\tau\notin \pbar(F)$,  
a {\it $\tau$-sequence}  is 
a sequence of distinct vertices 
$(v_1, v_2, \ldots, v_t)$  with $v_i\in \{s_{\beta+1},\ldots, s_q\}$  such that $\varphi(rv_1) =\tau$, and  the following three conditions are satisfied. 
\begin{enumerate} [(i)]
		\item $\{v_1,  \ldots, v_{t-1}\}$ is elementary and $\pbar(v_i) \notin \pbar(F)$ for each $i\in [1, t-1]$; 
		
		\item $\varphi(rv_i)=\pbar(v_{i-1})$ for each $i\in [2,t]$;   and 
		
		\item  There are three possibilities for  $\pbar(v_t)$:  (A)  $\pbar(v_t)=\tau$,  (B)  $\pbar(v_t )\in  \pbar(F)$, or (C) $\pbar(v_t) =\pbar(v_{i-1})$  for some $i\in [2,t-1]$. Accordingly, we name the $\tau$-sequence type A, type B, and type C, respectively, where a type A sequence is also called a {\em rotation}. 
	\end{enumerate} 

\end{DEF}	

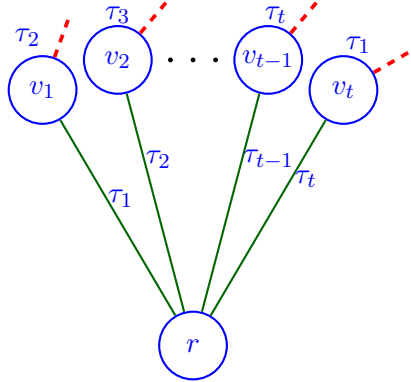
\begin{figure}[!htb]
	\begin{center}
		\begin{tikzpicture}[scale=1]
		
		{\tikzstyle{every node}=[draw ,circle,fill=white, minimum size=0.9cm,
			inner sep=0pt]
			\draw[blue,thick] (0, -3) node (r)  {$r$};
			\draw[blue,thick] (-2, 0.4) node (sa)  {$v_1$};
			\draw [blue,thick](-1, 0.8) node (sa2)  {$v_2$};
			\draw [blue,thick](1, 0.8) node (sb)  {$v_{t-1}$};
			\draw [blue,thick](2, 0.4) node (sb2)  {$v_t$};
		}
		\path[draw,thick,black!60!green]
		(r) edge node[name=la,pos=0.6] {\color{blue}\quad$\tau_1$} (sa)
		(r) edge node[name=la,pos=0.7] {\color{blue}\quad$\tau_2$} (sa2)
		(r) edge node[name=la,pos=0.7] {\color{blue}\quad\,\,\,\,\,$\tau_{t-1}$} (sb)
		(r) edge node[name=la,pos=0.7] {\color{blue}\quad$\tau_t$} (sb2);

		\draw[dashed, red, line width=0.5mm] (sa)--++(70:1cm); 
		\draw[dashed, red, line width=0.5mm] (sa2)--++(50:1cm); 
		\draw[dashed, red, line width=0.5mm] (sb)--++(50:1cm); 
		\draw[dashed, red, line width=0.5mm] (sb2)--++(30:1cm);

		\draw[blue] (-2.2, 1.1) node {$\tau_2$};  
		\draw[blue] (-1.0, 1.4) node {$\tau_3$};  
		\draw[blue] (1.1, 1.4) node {$\tau_{t}$}; 
		\draw[blue] (2.2, 1.0) node {$\tau_{1}$}; 
		
		{\tikzstyle{every node}=[draw ,circle,fill=black, minimum size=0.05cm,
			inner sep=0pt]
			
			\draw(-0.3,0.8) node (f1)  {};
			\draw(0,0.8) node (f1)  {};
			\draw(0.3,0.8) node (f1)  {};
			
		} 
		\end{tikzpicture}
	\end{center}
	\caption{A rotation in the neighborhood of $r$.}
	\label{rt}
\end{figure}

An example of a rotation is given in Figure~\ref{rt}, where $\tau_i = \pbar(v_{i-1})$ for
each $i\in [2, t]$.

\begin{LEM}\label{lem:tau-seq}
If $F$ is    maximum, then for any color $\tau\notin \pbar(F)$, there is a unique $\tau$-sequence.   
\end{LEM}

\pf Since $\tau\notin \pbar(r)$, there is a vertex $s\in N(r)$ such that $\varphi(rs) = \tau$. 
Since  $F$ is maximum, we have 
$\alpha+1,\beta+1\in 
\{\varphi(ru_1), \varphi(ru_2)\}$, and so $s \notin \{u_1, u_2\}$.  Since $\varphi(rs_i) \in \pbar(F)$ for all $i\in [2, \beta]$, 
$s =v_1$ for some $v_1\in \{s_{\beta+1}, \ldots, s_q\}$, where we recall $q=d(r)-2$. 

 Starting with a singlton sequence $(v_1)$, 
let $(v_1, \ldots, v_{t-1})$ be a longest sequence of vertices in $N(r)\setminus V(F)$ satisfying the following   two conditions:
\begin{enumerate} [(i)]
		\item $\{v_1,  \ldots, v_{t-1}\}$ is elementary and $\pbar(v_i) \notin \pbar(F)$ for each $i\in [1, t-1]$; and 
		
		\item $\varphi(rv_i)=\pbar(v_{i-1})$ for each $i\in [2,t-1]$. 
		
\end{enumerate}	
Let $v_t$ be a vertex in $N(r)$ such that $\varphi(rv_t) =\pbar(v_{t-1})$. 
Since $\pbar(v_{t-1})\notin \pbar(F)$, $v_t\in \{s_{\beta+1}, \ldots, s_q\}$. 
 If $v_t=v_1$,
then $(v_1, \ldots, v_{t-1})$ is a $\tau$-sequence of type A. Thus we assume that $v_t \ne v_1$, i.e.,
$\pbar(v_{t-1}) \ne \tau$. 
Since $\pbar(v_{t-1})\ne 	\pbar(v_i)$ for all $i\in [1,t-2]$, $v_t \notin \{v_2, \dots, v_{t-1}\}$. 
Hence $v_t\ne v_i$ for each $i\in [1,t-1]$. 
By the maximality of $(v_1, \dots, v_{t-1})$,  $\pbar(v_t)$ can only have three possibilities,  (A), (B) or (C),   as listed in condition (iii) of Definition 2.7. 

Moreover, since each $|\pbar(s_i)| =1$ for all $i\in [\beta+1, q]$, the  sequence above is unique. 
\qed

\begin{LEM}\label{lem:rs1}
If  $F$ is  maximum, then for any color  $\tau \notin \pbar(F)$, $r\in P_{s_1}(\tau, \Delta)$ and $r\in P_{s_1}(2,\tau)$. 
\end{LEM}


\pf We only show $r\in P_{s_1}(\tau, \Delta)$ since the proof of the other case is symmetric.  Assume to the contrary that $r\not\in P_{s_1}(\tau, \Delta)$. Let  $\varphi(rv_1)=\tau$ for $v_1\in \{s_{\beta+1}, \ldots, s_q\}$  and    $(v_1, \ldots, v_t)$  
be the $\tau$-sequence by Lemma~\ref{lem:tau-seq}. 
Let 
$\varphi' = \varphi/C_r(\tau, \D)$. Notice that under the coloring $\varphi'$, 
$\varphi'(rv_1) =\D$ and $\varphi'(rs_{\alpha +1}) = \tau$,  and the color on each edge from $E(F)\setminus\{rs_{\alpha+1}\}$ and the missing color on each vertex of $F$ are the same as the corresponding colors under $\varphi$. Hence, $F'  = (r,rs_1,s_1,\ldots, rs_\alpha,s_\alpha, rv_1, v_1, \ldots, rv_t, v_t)$ is a multifan w.r.t. $e$ and $\varphi'$.  
 
    We consider three cases according to the type of this $\tau$-sequence: type A, type B, or type C with respect to the coloring $\varphi$.   
We note that if $\pbar'(v_t) =\tau$, then, due to $\varphi'(rs_{\alpha+1}) = \tau$,   $F'$ can be extended to a larger multifan:
 $$
 F^* = (r,rs_1,s_1,\ldots, rs_\alpha,s_\alpha, rv_1, v_1, \ldots, rv_t, v_t, rs_{\alpha+1}, s_{\alpha+1}, \ldots, rs_\beta, s_\beta)$$
 which is also larger than $F$, giving a contradiction to the maximality of $F$.  We will use $F'$ and $F^*$ to lead a contradiction in our proof.

{\flushleft \bf Type A:}   In this case $\pbar(v_t)=\tau$.    If $C_r(\tau, \D, \varphi) \ne P_{v_t}(\tau,\Delta,\varphi)$, then 
$ \pbar'(v_t)= \tau$, and so $F^*$ is a multifan w.r.t. $e$ and $\varphi'$,  giving a contradiction.   Thus $C_r(\tau, \D, \varphi) =P_{v_t}(\tau,\Delta,\varphi)$, which in turn gives $\pbar'(v_t) = \D$. 
In this case, $\Delta \in \pbar'(s_1)\cap \pbar'(v_t)$, and so $F'$ is not elementary, giving a contradiction.

{\flushleft \bf Type B:}  In this case $\pbar(v_t)$, denoted by $\gamma$, is in $\pbar(F)$. 
If $\gamma \in \pbar(\{r, s_1, \dots, s_{\alpha}\})$,  then $F'$ is not elementary, giving a contradition. 
Thus, we have either $\gamma =\D$  or $\gamma \in \pbar(\{s_{\alpha +1}, \dots, s_{\beta})$. 

Assume first that $\gamma \ne \Delta$. Let $\gamma= \pbar(s_j)$ for some $j\in [\alpha+1, \beta]$. 
Since $P_r(1,\gamma,\varphi)=P_{s_j}(1,\gamma, \varphi)$ 
and $1,\gamma  \not\in \{\tau,\Delta\}$, we still have $P_r(1,\gamma ,\varphi') = P_r(1, \gamma, \varphi) = P_{s_j}(1,\gamma ,\varphi')$. 
Let  $\varphi''=\varphi'/P_{v_t}(1,\gamma , \varphi')$.  
Under $\varphi''$,  $F'$ is also  a multifan.
However, color $1\in \pbar''(r)\cap \pbar''(v_t)$, giving a contradiction to  $V(F')$ being elementary.  
Thus $\gamma =\Delta$. If $C_{r}(\tau,\Delta,\varphi)=P_{v_t}(\tau,\Delta,\varphi)$, then $\pbar'(v_t) = \tau$, which in turn shows that $F^*$ is a multifan w.r.t. $e$ and $\varphi'$, a contradiction. 
Thus $C_{r}(\tau,\Delta,\varphi)\ne P_{v_t}(\tau,\Delta,\varphi)$. So $\varphi'(v_t) =\D$, which in turn shows that $F'$ is not elementary since $\D$ is also in $\pbar'(s_1)$, a contradiction.

{\flushleft \bf Type C:}  Suppose  $\pbar(v_t) =\pbar(v_{i-1})=\tau_i$  for some $i\in [2,t-1]$ and some $\tau_i\in [1,\Delta]\setminus \pbar(F)$.
Note that one of $v_{i-1}$ and $v_{t}$ is  $(1,\tau_i)$-unlinked with $r$.
By doing a $(1,\tau_i)$-swap at a vertex in $\{v_{i-1}, v_t\}$
that is $(1,\tau_i)$-unlinked with $r$, 
we convert this case to the Type B case. 
\qed

Note that under the condition of Lemma~\ref{lem:rs1}, if $P$ is a $(2,\tau)$- or $(\tau, \Delta)$-chain disjoint from $P_{s_1}(2,\tau )$ or $P_{s_1}(\tau, \D)$, then we also have $r\notin P$, and so 
the Kempe change $\phiv/P$ gives an $F$-stable coloring.

\begin{DEF}\label{def-shift}
	Let $h,\ell \in [1,q]$. 
The \emph{shifting from $s_h$ to $s_\ell$}  is a recoloring  operation $rs_{i}:\  \varphi(rs_{i}) \rightarrow \pbar(s_{i})$ for all
 $ i\in[h,\ell]$, i.e., replacing the current color on the edge $rs_{i}$ with the missing color at $s_{i}$ for all $i\in [h, \ell]$. 
 \end{DEF}
 
 We apply shiftings when the sequence $(s_h,\ldots, s_\ell)$ forms a rotation or 
 is a B-type $\tau$-sequence, where $\tau=\varphi(rs_h)$, such that $\pbar(s_\ell)=1$. 
 Since $1\in \pbar(r)$, we obtain another $\D$-edge coloring in both cases. 
 We do not know whether a shifting can be achieved through a sequence of Kempe changes. So, in this paper, ``Kempe changes" do not include ``shifting".    In the proof, we sometimes use the following weaker version of ``stable" coloring. 
 
 \begin{DEF}\label{defi:wstable}
A coloring $\phiv'\in \CC^\Delta(G-rs_1)$ is {\em $V(F -r)$-stable}  (w.r.t.  $F$ and $\phiv$) if 
$V(F)$ is the vertex set of a multifan $F_{\phiv'}$ at $r$ w.r.t. $e=rs_1$ and $\phiv'$,   
$\pbar'(s_1)=\pbar(s_1)=\{2,\Delta\}$, and $\pbar'(V(F_{\phiv'})\setminus \{r\})=\pbar(V(F)\setminus \{r\})$. 
 Moreover, a $V(F-r)$-stable coloring $\phiv'$ is called {\em  $V(F)$-stable} if $\pbar'(r)=\pbar(r)$.
 \end{DEF} 
 
 \begin{LEM}\label{Lem:Wstable}
 For any color $\gamma\in \pbar(F)$ and a vertex $x\in V(G)\setminus V(F)$, the following two statements hold. 
 \begin{itemize}
 \item if $\gamma$ is $2$-inducing, then the Kempe change $\phiv/P_x(\gamma, \D, \phiv)$ gives a $V(F)$-stable coloring provided $\pbar(x)\cap \{\gamma, \D\} \ne \emptyset$, and 
 \item if $\gamma$ is  $\D$-inducing, then the Kempe change $\phiv/P_x(2,\gamma, \phiv)$ gives a $V(F)$-stable coloring provided $\pbar(x)\cap \{\gamma, 2\} \ne \emptyset$. 
 \end{itemize}
 \end{LEM}
 
 \pf By symmetry, we only prove the first one. If $r\notin P_x(\gamma, \D, \phiv)$, we are done by Lemma~\ref{Lem:Stable}.  Assume $r\in P_x(\gamma, \D, \phiv)$.  Since $\pbar(x)\cap\{\gamma, \D\}\ne \emptyset$, $P_x(\gamma, \D, \phiv)$ is disjoint from $P_{s_1}(\gamma, \D, \phiv) = P_{\pbar^{-1}_F}(\gamma, \D, \phiv)$. Let $\phiv' = \phiv/P_x(\gamma, \D, \phiv)$. Note that $\varphi (rs_{\alpha+1}) = \D$.  Let $s_i=\pbar^{-1}_F(\gamma) $ for some $i\in [1,\alpha]$.  We have $\phiv'(rs_{\alpha +1}) =\gamma$ and $\phiv'(rs_i) =\D$. So, 
 $$
 F' = (r, rs_1, s_1, rs_2, \dots, rs_{i-1}, rs_{\alpha +1}, s_{\alpha +1}, \dots, s_\beta, rs_i, \dots, s_{\alpha} )
 $$ 
 is a multifan w.r.t. 
 $rs_1$ and $\phiv'$. Clearly, $\pbar'(s_1) =\pbar(s_1) = \{2, \D\}$ and $\pbar'(V(F_{\phiv'})\setminus \{r\})=\pbar(V(F)\setminus \{r\})$. Hence, $\phiv'$ is $V(F)$-stable. \qed

 Let $\tau\in [1,\D]\setminus \pbar(F)$ and $(v_1, v_2, \ldots, v_t)$ 
 be the $\tau$-sequence at $r$. 
If the $\tau$-sequence  
is of type A,  the shifting of this sequence yields a coloring 
in $\CC^{\D}(G-e)$, which is   $F$-stable.  We call such an operation an {\it $A$-shifting}.   If the $\tau$-sequence  is of  type B and satisfies  $\pbar(v_t) = 1$,  the shifting of this  sequence yields a coloring 
$\varphi'\in \CC^{\D}(G-e)$ with $\pbar'(r) = \tau$, which is $V(F-r)$-stable.  We call such an operation a {\it $B$-shifting}.

Let $P$ be a $(\tau, *)$-chain with endvertices $x$ and $y$,
where $*$ represents any color from $[1,\Delta]\setminus \{\tau\}$. 
  Suppose that  $rv_1 \in E(P)$ and  $x, y\notin \{v_1, \dots, v_t\}$. If either the $A$-shifting or the $B$-shifting is eligible, we do it and obtain a new coloring $\phiv'$. Notice that $\pbar'(v_1) =\tau$. So, either $P_x(\tau, *, \phiv') = P_{v_1}(\tau, *, \phiv')$ or $P_y(\tau, *, \phiv') = P_{v_1}(\tau, *, \phiv')$ but not both. Consequently, $x$ and $y$ are $(\tau, *)$-unlinked w.r.t. coloring $\phiv'$. We will use this ``unlink" technique in the following lemma. In the proofs, we may need to preserve some colors at
  a vertex, which leads to the following definition.

\begin{DEF}\label{def-avoid}
Given a set $S$ of colors, a coloring $\varphi'$ is called {\em $S$-avoiding} (w.r.t. $\varphi$) if every Kempe change  applied in obtaining $\varphi'$ from coloring $\varphi$ does not involve any color from $S$. 
\end{DEF}

In the following lemma,  whenever $P_x(\tau, \D)$
or $P_x(2,\tau)$ is used, it 
 implicitly  implies that  one of the two colors from $\tau$ and $\D$ or from  $2$
 and $\tau$ is missing at $x$. 
\begin{LEM}\label{lem:tau-seq2}
Suppose $F$ is a maximum multifan,  $N[r] \ne V(G)$ and $\pbar(V(F)) \ne [1, \D]$.  For any  vertex   $x\in V(G)\setminus N[r]$ and any  color $\tau \in [1,\D]\setminus \pbar(F)$, let $v_t$ be the last vertex of the $\tau$-sequence.   If   $\pbar(x) \cap \{\tau, \D\} \ne \emptyset$,  then 
\begin{enumerate}[(i)]
	\item  If  $\tau\in \pbar(x)$, then  there is an $F$-stable coloring $\varphi'\in \CC^\Delta(G-rs_1)$  such that $1\in \pbar'(x)$.

\item  If $\tau\in \pbar(x)$, then there exists a $V(F-r)$-stable and $\{\D\}$-avoiding coloring $\phiv'$ such that 
$\pbar'(r) \in \pbar'(x)$ unless  the $\tau$-sequence is of type $B$ and  $\pbar(v_t) =\D$. 


	
\item There exists a $V(F -r)$-stable and  $\{\tau,\Delta\}$-avoiding coloring $\varphi'$ such that 
$P_{s_1}(\tau,\Delta, \phiv')\ne P_{x}(\tau,\Delta, \phiv')$  unless the $\tau$-sequence is of type B and  $\pbar(v_t)=\Delta$.

\item There exists a $V(F -r)$-stable  and  $\{2, \tau,\Delta\}$-avoiding coloring $\varphi'$ such that 
$P_{s_1}(\tau,\Delta, \phiv')\ne P_{x}(\tau,\Delta, \phiv')$  unless the $\tau$-sequence is of type B and  $\pbar(v_t)\in \{2, \Delta\}$. 
	   
\item There exists a  $V(F)$-stable and $\{1,\tau\}$-avoiding coloring $\varphi'$ such that 
$P_{s_1}(\tau,\Delta, \phiv')\ne P_{x}(\tau,\Delta, \phiv')$ or  $P_{s_1}(2,\tau, \phiv')\ne P_{x}(2,\tau, \phiv')$  unless 
the $\tau$-sequence is of type B and $\pbar(v_t)=1 $.

\item  There exists a $V(F)$-stable and  $\{1, \tau, \Delta\}$-avoiding $\varphi'$ such that 
$P_{s_1}(\tau,\Delta, \phiv')\ne P_{x}(\tau,\Delta, \phiv')$ unless the $\tau$-sequence is of type B and  either $\pbar(v_t) \in \{1, \D\}$   or   $\pbar(v_t)$ is 2-inducing. 
	  		
\item If $\tau\in \pbar(x)$, then there exists a $V(F)$-stable and  $\{1,\tau, \Delta\}$-avoiding coloring $\varphi'$ such that 
$P_{s_1}(2,\tau, \phiv')\ne P_{x}(2,\tau, \phiv')$  unless the $\tau$-sequence is of type B and  either $\pbar(v_t) \in \{1, \D\}$   or    $\pbar(v_t)$ is 2-inducing.

\end{enumerate}
\end{LEM}

\pf  
Let $(v_1, \ldots, v_t)$ be the $\tau$-sequence at $r$.  
By Lemma~\ref{lem:rs1}, $r\in P_{s_1}(\tau,\Delta,\varphi)$
and $r\in P_{s_1}(2,\tau,\varphi)$. For the proof of each statement in (iii)-(vii),
we may assume that $ P_x(\tau,\Delta,\varphi)=P_{s_1}(\tau,\Delta,\varphi)$
or $P_x(2,\tau,\varphi)= P_{s_1}(2,\tau,\varphi)$ in accordance with its assertion. 
We will apply either an A-shifting or a B-shifting on $(v_1, \ldots, v_t)$,
which will cutoff the linkage between $x$ and $s_1$ by the remark prior to Definition~\ref{def-avoid}.

Assume first that $(v_1, \ldots, v_t)$ is of type A. For each of (iii) to (vi), since $F$ is maximum, $r\in P_{s_1}(\tau, \D)$, so under coloring $\phiv$,  edge $rv_1\in P_{s_1}(\tau, \D) = P_x(\tau, \D)$. 
We do the $A$-shifting on the $\tau$-sequence to obtain an $F$-stable coloring $\phiv'$. Clearly $\phiv'$ is $V(F)$-stable  and  $P_{s_1}(\tau, \D,\phiv') \ne P_{x}(\tau, \D, \phiv')$, as desired.  For (vii), the argument is the same as above with the color 2 playing the role of $\D$.  Now consider (i) and (ii). For (i), we assume instead that $x$ and $r$ are $(1,\tau)$-linked under $\phiv$, since otherwise we just do a $(1,\tau)$-swap at $x$ to get a desired coloring. Similarly, for (ii), we assume that $x$ and $r$ are $(1,\tau)$-linked under $\phiv$. For both (i) and (ii), we do a $(1,\tau)$-swap at $v_t$, and then do the shifting from $v_1$ to $v_t$ to obtain coloring $\phiv'$. Note that $\phiv'$ is $V(F-r)$-stable, and $\pbar'(r)=\tau\in \pbar'(x)$ now. So $\phiv'$ is a desired coloring of (ii). For (i), we obtain $\phiv''$ from $\phiv'$ by renaming $\tau$ as $1$ and vice versa. Then $\phiv''$ is $F$-stable and it is a desired coloring of (i).

Assume now that $(v_1, \ldots, v_t)$ is of type B. Let $\pbar(v_t)=\gamma\in \pbar(F)$. Recall that $r$ and $\pbar^{-1}_F(\gamma)$ are $(1,\gamma)$-linked by Lemma~\ref{thm:vizing-fan2}~\eqref{thm:vizing-fan1b}.  For (i),  we simply do a $(1,\gamma)$-swap at $v_t$ and then do a shifting from $v_1$ to $v_t$. Then by renaming $\tau$ as $1$ and vice versa, we obtain an $F$-stable coloring $\varphi'$ such that color $1$ is missing at $x$.  For (ii), we are done if $\gamma=\Delta$. Otherwise we do the same as in the proof for (i) without renaming $\tau$ as $1$. The resulting coloring is $V(F-r)$-stable and is $\{\Delta\}$-avoiding,  and  $\pbar'(r)=\tau\in \pbar'(x)$. For (iii), if  $\gamma  \ne \Delta$, we still do a $(1,\gamma)$-swap at $v_t$ and then do a shifting from $v_1$ to $v_t$. The resulting coloring is $V(F -r)$-stable and is $\{\tau,\Delta\}$-avoiding,  and  $P_{s_1}(\tau,\Delta)\ne P_{x}(\tau,\Delta)$ under the current coloring.  We do the same argument for (iv) except that $\gamma \in \{2,\Delta\}$.  For (v),  we may assume that $\gamma \ne 1$ and $\gamma$ is $2$-inducing by the symmetry between 2 and $\Delta$. Since $s_1$ and $\pbar_{F}^{-1}(\gamma)$ are $(\gamma,\Delta)$-linked by Lemma~\ref{thm:vizing-fan2}~\eqref{thm:vizing-fan2-a}, we first do a $(\gamma, \Delta)$-swap at $v_t$. The resulting coloring is $V(F)$-stable, so we still have $r\in P_{s_1}(\tau,\Delta)=P_{x}(\tau,\Delta)$. Now we do a $(\tau,\Delta)$-swap at $v_t$, this results in a type A $\tau$-sequence. For (vi) and (vii), we are done if $\gamma$ is $1$ or $\Delta$ or 2-inducing.  If $\gamma$ is $\Delta$-inducing, we do a $(2,\gamma)$-swap at $v_t$. The resulting coloring is $V(F)$-stable by Lemma~\ref{thm:vizing-fan2}~\eqref{thm:vizing-fan2-a}, and it is $\{1,\tau,\Delta\}$-avoiding. Now the missing color of $v_t$ is a 2-inducing color, as desired.

Assume finally that $(v_1, \ldots, v_t)$ is of type C. That is, $\pbar(v_t)= \pbar(v_{i-1})=\tau_i$  for some $i\in [2,t-1]$ and some $\tau_i\in [1,\D]\setminus \pbar(F)$.   We show that under each assumption of the statements, we can reduce this sequence into a B-type $\tau$-sequence. For each of (i) to (iv), since one of $v_{i-1}$ and $v_t$ is $(1,\tau_i)$-unlinked with $r$, we do a $(1,\tau_i)$ swap at a vertex in $\{v_{i-1}, v_t\}$ that is $(1,\tau_i)$-unlinked with $r$, resulting in a type B $\tau$-sequence $(v_1, \ldots, v_{i-1})$ or  $(v_1, \ldots, v_t)$. For (v), since $\tau_i\in [1,\Delta]\setminus \pbar(F)$, by Lemma~\ref{lem:rs1}, $r\in P_{s_1}(\tau_i,\Delta)$. Since  one of  $v_{i-1}$ and $v_t$ is $(\tau_i, \Delta)$-unlinked with $r$, we do a $(\tau_i,\Delta)$-swap at a vertex in $\{v_{i-1}, v_t\}$ that is $(\tau_i, \Delta)$-unlinked with $r$, resulting in a type B $\tau$-sequence $(v_1, \ldots, v_{i-1})$ or  $(v_1, \ldots, v_t)$. For (vi) and (vii), since $\tau_i\in [1,\Delta]\setminus \pbar(F)$, by Lemma~\ref{lem:rs1}, $r\in P_{s_1}(2,\tau_i)$. Since one of  $v_{i-1}$ and $v_t$ is $(2,\tau_i)$-unlinked with $r$, we do a $(2,\tau_i)$-swap at a vertex in $\{v_{i-1}, v_t\}$ that is $(2, \tau_i)$-unlinked with $r$, resulting in a type B $\tau$-sequence $(v_1, \ldots, v_{i-1})$ or  $(v_1, \ldots, v_t)$.\qed

\section{Proof of Theorem~\ref{thm:s1-adj}}

\begin{THM2}\label{thm:s1-adj1}
	Let $G$ be a  class 2 graph with maximum degree $\Delta$, 
	$r\in V_{\Delta}$ be light, and $s\in N(r)$ with $d(s)<\Delta$.  If $rs$ is a critical edge 
	of $G$, then all vertices in $N(s) \setminus N(r)$ are $\D$-vertices. 
\end{THM2}

\pf Assume to the contrary that there exists $x\in N(s)\setminus N(r)$ with $d(x) < \D$.
Clearly, $x\ne r$.  
Denote  $s$ by $s_1$.  
Let $\varphi\in \CC^\Delta(G-rs_1)$  and assume the corresponding multifan  $F$ w.r.t. $rs_1$  is maximum and typical.

We claim that there is an $F$-stable coloring such that color 1 is missing at $x$. To see this, 
let $\tau \in \pbar(x)$. If $\tau\in \pbar(F)$, then $\pbar_F^{-1}(\tau)$ and $r$
are $(1,\tau)$-linked by Lemma~\ref{thm:vizing-fan2}~\eqref{thm:vizing-fan1b}. So, $P_x(1, \tau)$ does not contain any edge of $F$ and does not end at any vertex in $F$. Hence $\phiv/P_x(1, \tau)$ is $F$-stable such that color $1$ is missing at $x$. 
 We assume that $\tau \notin  \pbar(F)$. 
By Lemma~\ref{lem:tau-seq2} (i), there is an $F$-stable coloring such that color 1 is missing at $x$.   So we assume  $1\in \pbar(x)$.

Let $\varphi(s_1x)=\tau$.  If $\tau\in \pbar(F)$, 
we may assume it is 2-inducing.  Since $\pbar_F^{-1}(\tau)$ and $r$
are $(1,\tau)$-linked by Lemma~\ref{thm:vizing-fan2}~\eqref{thm:vizing-fan1b}, 
we  do a $(1,\tau)$-swap at $x$ and  get an $F$-stable coloring.  We then do a $(\tau,\Delta)$-swap at $x$ and get a new coloring $\phiv'$. 
 Since $\tau$ is $2$-inducing, $s_1$ and $\pbar_F^{-1}(\tau)$ are $(\tau, \D)$-linked,   $\D$ is still missing at $s_1$. 
We see that $F^*=(r,rs_1,s_1, s_1x,x)$ is a multifan w.r.t. $\phiv'$.  However,
we have $ \Delta \in \pbar'(s_1)\cap \pbar'(x)$, contradicting  $V(F^*)$ being elementary. Thus we assume that $\tau \notin  \pbar(F)$.
We do a $(1,\Delta)$-swap at $x$ and get an $F$-stable coloring $\phiv'$. 
Then $P_{s_1}(\tau,\Delta)=s_1x$ that does not contain vertex $r$,  showing a contradiction to Lemma~\ref{lem:rs1}. 
\qed

\section{Proof of Theorem~\ref{thm:longk}}
 In this section, we assume that  $G$ is a  class 2 graph with maximum degree $\Delta$, $r\in V(G)$ is a {\it light} vertex and 
 $s_1\in N_{\D -1}(r)$ such that $rs_1$ is a critical edge, and vertex $x\in V(G)\setminus N[r]$ with  $d(x) \le \D -2$.

 \begin{THM3}\label{thm:longka}
	If $d(r) =\D$ and  $d(x) \le \D -3$, then $N(x)\cap N(s_1) \subseteq N(r)\setminus N_\Delta(r)$. 
\end{THM3}

The proof of Theorem~\ref{thm:longk} is based on the following three lemmas whose proofs will be given in the following three subsections, respectively.  Let  $\phiv\in \CC^{\D}(G-rs_1)$ and $F$ be a typical multifan w.r.t. $rs_1$ and $\phiv$.   We additionally assume that $F$ is a maximum multifan w.r.t. edge $rs_1$. 

\begin{LEM}\label{lem:good-coloring0}
Suppose $r$ is a $\D$-vertex.  For every vertex $u\in N(s_1)\cap N(x)$ such that $u\not\in N(r)\setminus N_{\Delta}(r)$,
 if $\{1, 2\}\subseteq \pbar(x)$ then $\phiv(ux) \ne \D$ and if $\{1, \D\} \subseteq \pbar(x)$ then $\phiv(ux) \ne 2$.   
\end{LEM}

\begin{LEM}\label{lem:good-coloring}
Suppose $r$ is a $\D$-vertex. For every vertex $u\in N(s_1)\cap N(x)$ such that $u\not\in N(r)\setminus N_{\Delta}(r)$, if $\{2,\Delta\} \subseteq \pbar(x)$ then  $\varphi(s_1u)\ne 1$. 
 \end{LEM}	

By relaxing the condition $d(r) =\D$ to $d(r) \ge \D -1$, we have the following results. 
\begin{LEM}\label{lem:good-coloring1}
Under the assumption  $d(r) \ge \D -1$, the following statements hold. 
\begin{enumerate}[(i)]
\item If $d(x)\le \Delta-3$ and $1 \in \pbar(x)$,  then there is a $V(F)$-stable coloring 
		such that  both colors $2$ and $\D$ are missing at $x$.
		
\item If  $\{2,\Delta\} \subseteq  \pbar(x)$, then 
for every vertex $u\in N(s_1)\cap N(x)$ such that $u\not\in N(r)\setminus N_{\Delta}(r)$, 	we have the following two claims: 

\begin{enumerate}[(a)]	
\item there is a $V(F)$-stable coloring  $\varphi'$ 
such that   $1\in \pbar'(x)$ and $\varphi'(ux)\in \{2,\Delta\}$,  and, more specifically,  $\varphi'(ux)=\Delta$ if $\varphi'(s_1u)$ is 2-inducing and  $\varphi'(ux)=2$  if $\varphi'(s_1u)$ is $\Delta$-inducing;  and

\item    there is  a $V(F)$-stable coloring  $\phiv'$ 	such that 
$\varphi'(s_1u)=1$ and 
$\{2, \Delta\}\cap  \pbar'(x) \ne \emptyset$. Moreover, $\{2, \D\}\subseteq \pbar'(x)$ if additionally  $d(x) \le \D -3$. 
\end{enumerate}

\end{enumerate}
\end{LEM}

{\flushleft \bf Proof of Theorem~\ref{thm:longk}.} 
Let $N_{\Delta-1}(r)=\{s_1, s_2,\ldots, s_{\Delta-2}\}$. 
 Suppose  to the contrary that there is a vertex $u\in V(G)$ such that $u\not\in N(r)\setminus N_\Delta(r)$ and $u$ is adjacent to both $x$ and $s_1$.  Then $u\not\in \{r,s_1,\ldots, s_{\Delta-2}\}$ since $x\not\in N(r)$.
Following the notation given at the beginning of this section,  we let $\varphi\in \CC^\Delta(G-rs_1)$, and $F$ be a typical multifan w.r.t. $rs_1$ and $\phiv$. We also assume that $F$ is a maximum multifan w.r.t. $rs_1$.

We claim that there exists a $V(F)$-stable coloring $\phiv'$ such that $1 \in \pbar'(x)$.  Let $\tau \in \pbar(x)$.  If $\tau\in \pbar(F)$,  then by Lemma~\ref{Lem:Stable} the Kempe change $\phiv/P_x(1, \tau, \phiv)$ gives  an $F$-stable coloring $\phiv'$. Clearly,  $1\in \pbar'(x)$.  Thus, we assume that $\tau \in [1,\Delta]\setminus \pbar(F)$. 
By Lemma~\ref{lem:tau-seq2} (i), there is an $F$-stable coloring $\phiv'$ such that $1\in \pbar'(x)$. 

Applying  Lemma~\ref{lem:good-coloring1} (i), and then (ii) with $\phiv'$,  we get a $V(F)$-stable coloring $\phiv''$ such that 
$\{ 2, \D \}\subseteq \pbar''(x)$ and $\phiv''(s_1u) =1$, which gives a contradiction to Lemma~\ref{lem:good-coloring}. 
\qed

\subsection{Proof of Lemma~\ref{lem:good-coloring0} }

 By symmetry, we only prove the first conclusion.  
 Suppose to the contrary  and,  without loss of generality,   that   $\{1,2\} \subseteq \pbar(x)$ and 
$\varphi(ux)=\Delta$.  Notice that $1 \in \pbar(r)\cap \pbar(x)$ and $\D \in \pbar(s_1)\cap \pbar(x)$.  Let $\tau=\varphi(s_1u)$. 

Consider first that $\tau \in \pbar(F)$.   If $\tau = 1$, then $\phiv(s_1u) \in \pbar(r)$ and $\phiv(ux) = \D \in \pbar(s_1)$, and so 
$K=(r,rs_1,s_1,s_1u, u, ux, x)$ is a Kierstead path. Since $d(s_1) =\D -1$, $V(K)$ is elementary by Lemma~\ref{Lemma:kierstead path1}, showing a contradiction to
$2\in \pbar(x)\cap \pbar(s_1)$.  So, $\tau \ne 1$. Let $\varphi^*=\varphi/P_x(1,\tau, \phiv)$.  By Lemma~\ref{Lem:Stable}, $\phiv^*$ is $F$-stable. 
Suppose that $\tau$ is 2-inducing. 
 If $s_1u\not\in P_x(1,\tau,\varphi)$, then $\phiv^*(s_1u) = \tau$, and so $P_{s_1}(\tau, \Delta, \phiv^*)=s_1ux$,
contradicting that $s_1$ and $\pbar^{-1}_F(\tau)$ are $(\tau, \D)$-linked (Lemma~\ref{thm:vizing-fan2}~\eqref{thm:vizing-fan2-a}). 
If  $s_1u\in P_x(1,\gamma,\varphi)$, then $\phiv^*(s_1u) =1\in \pbar^*(r)$ and $\phiv^*(ux) =\D\in \pbar^*(s_1)$, so under $\phiv^*$, 
$K=(r, rs_1,s_1,s_1u, u, ux,x)$  is a Kierstead path with $d(s_1) =\D -1 < \D$, but 2 is missing at both $s_1$
and $x$, showing a contradiction to $V(K)$ being elementary (Lemma~\ref{Lemma:kierstead path1}). 
Thus,  $\tau$ is $\Delta$-inducing. We do  $(\Delta,1)-(1,2)$-swaps at $x$ and get an $F$-stable coloring $\phiv'$ (Lemma~\ref{Lem:Stable}).   
Notice that $\phiv'(ux) =2$ and $\{1, \D\}\subseteq  \pbar'(x)$. 
This gives back to the 
previous case by the symmetry between 2 and $\Delta$, which leads to a contradiction.  Thus, $\tau\in [1,\Delta]\setminus \pbar(F)$.

Since $F$ is a maximum multifan, by Lemma~\ref{lem:tau-seq}
there is a unique $\tau$-sequence  $(v_1, \ldots, v_t) $.  We claim  that 
$s_1u\in P_{x}(1,\tau, \phiv)=P_r(1,\tau, \phiv)$.  Otherwise, 
let $\phiv' = \phiv/P_x(1, \tau, \phiv)$. Clearly, $\tau \in \pbar'(x)$. If $s_1u \notin P_x(1, \tau, \phiv)$, then 
$P_{s_1}(\tau, \D,\phiv')=s_1 u x$.  In this case, if $P_x( \tau,\D, \phiv)$ does not end at $r$, then 
$\phiv'$ is $F$-stable, which in turn gives $r\in P_{s_1}( \tau,\D, \phiv') = s_1ux$ by Lemma~\ref{lem:rs1}, a contradiction; 
if  $P_x( \tau, \D,\phiv)$ ends at $r$, then $\phiv'$ is $V(F-r)$-stable and $\pbar'(r) = \tau$, which in turn gives 
$P_{s_1}( \tau,\D, \phiv') =s_1ux$, which should contain $r$ 
 and end at $r$  by Lemma~\ref{thm:vizing-fan2}~\eqref{thm:vizing-fan1b}, giving a contradiction. So,  $P_x(1, \tau, \phiv)$ contains edge $s_1u$ and does not end at $r$. In this case, $\phiv'(s_1u) =1 \in \pbar'(r)$ and $\phiv'(ux) = \phiv(ux) =\D \in \pbar'(s_1)$,  and so 
 $K'=(r, rs_1, s_1, s_1u, u, ux, x)$ is a Kierstead path. But, $2\in \pbar'(s_1)\cap \pbar'(x)$ shows that $V(K')$ is not elementary, a contradiction. 

We consider below the $\tau$-sequence $(v_1, \ldots, v_t) $ according to its type, but deal with the situation in the following claim first. 
\begin{CLA} \label{cla:noB1} There does not exist a $V(F)$-stable coloring $\phiv'$ with $\varphi'(s_1u)=\tau$, 
	$\phiv'(ux) = \D$, $2\in \pbar'(x)$, and 
	 the $\tau$-sequence w.r.t. it is of $B$-type with $\pbar'(v_t) =1$. 
\end{CLA} 

\pf Suppose to  the contrary that there is such a $V(F)$-stable coloring. We also assume that under coloring $\phiv'$,  the $\tau$-sequence is also $(v_{1}, \dots, v_{t})$.  We do the $B$-shifting from $v_1$ to $v_t$ to get a new coloring $\phiv^*$. Note that $\phiv^*$ is a $V(F-r)$-stable coloring,  and  $\phiv^*(s_1u) =\tau = \pbar^*(r) $ and $\phiv^*(ux)=\D \in \pbar^*(s_1)$, which in turn shows that  $K=(r, rs_1, s_1, s_1u, u, ux, x)$ is a Kierstead path. But, $2\in \pbar^*(s_1)\cap \pbar^*(x)$ shows that $V(K)$ is not elementary, a contradiction.  \qed

If the $\tau$-sequence is of A-type, i.e., $\pbar(v_t)=\tau$, 
we do a $(1,\tau)$-swap at 
$v_t$ to get a coloring $\phiv'$.  
Since $s_1u\in P_{x}(1,\tau, \phiv)=P_r(1,\tau, \phiv)$,   $\phiv'$ is  $F$-stable. We also notice that $\phiv'(ux) =\phiv(ux) =\D$ and $2 \in \pbar'(x)$, which gives a contradiction to Claim~\ref{cla:noB1}.

Suppose that  the $\tau$-sequence is of C-type,  more specifically,  $\pbar(v_t)= \pbar(v_{i-1})=\tau_i$  for some $i\in [2,t-1]$. Since  one of $v_{i-1}$ and $v_{t}$ is $(1,\tau_i)$-unlinked with $r$, we do a $(1,\tau_i)$-swap at a vertex in $\{v_{i-1}, v_t\}$
that is $(1,\tau_i)$-unlinked with $r$ to get an $F$-stable coloring $\phiv'$.  Clearly, $1\in \pbar'(v_{i-1})$ or $1\in \pbar'(v_t)$.  In either case, the resulting $\tau$-sequence is of type B with color 1 missing at the last vertex, which gives a contradiction to Claim~\ref{cla:noB1}. 

Suppose  now that $\tau$-sequence is of B-type and let   $\pbar(v_t) =\gamma$. Then, $\gamma\in \pbar(F)$.   By Claim~\ref{cla:noB1}, $\gamma \ne 1$.  
If $\gamma\ne \Delta$, we first do a $(1,\gamma)$-swap at $v_t$ and get an $F$-stable coloring $\phiv'$. Note that $1\notin \pbar'(x)$ may occur. Under coloring $\phiv'$, the $\tau$-sequence is of type B and $1\in \pbar'(v_t)$,
giving a contradiction to Claim~\ref{cla:noB1}.  
Thus, $\gamma=\Delta$. 
We consider two cases regarding whether $t=1$. 

\medskip 

{\bf \noindent Case 1. $t=1$}. 

\begin{figure}[!htb]
	\begin{center}
		\begin{tikzpicture}[scale=1]

	{\tikzstyle{every node}=[draw ,circle,fill=white, minimum size=0.5cm,
		inner sep=0pt]
		\draw[blue,thick](-3,-2) node (s1)  {$s_1$};
		\draw[blue,thick](-3,-4) node (u)  {$u$};
		\draw[blue,thick](-3,-6) node (x)  {$x$};
		
		\draw[blue,thick](-1.5,-4) node (r)  {$r$};
		\draw[blue,thick](-1,-2.5) node (v1)  {$v_1$};
	}
	\path[draw,thick,black!60!green]
	(s1) edge node[name=la,pos=0.6, above] {\color{blue} $\tau\quad$} (u)
	(u) edge node[name=la,pos=0.6, above] {\color{blue}$\Delta\quad$} (x)
	(r) edge node[name=la,pos=0.4, above] {\color{blue} $\tau\quad$} (v1);
	
	\draw[dashed, red, line width=0.5mm] (s1)--++(20:1cm); 
	\draw[dashed, red, line width=0.5mm] (s1)--++(160:1cm); 
	\draw[dashed, red, line width=0.5mm] (x)--++(340:1cm); 
	\draw[dashed, red, line width=0.5mm] (x)--++(200:1cm);
	\draw[dashed, red, line width=0.5mm] (r)--++(160:1cm);
	\draw[dashed, red, line width=0.5mm] (v1)--++(20:1cm); 
	
	\draw[blue] (-3.5, -2.1) node {$2$};
	\draw[blue] (-2.5, -2.1) node {$\D$};  
	\draw[blue] (-2.1, -4) node {$1$};  
	\draw[blue] (-0.6, -2.1) node {$\D$}; 
	\draw[blue] (-3.5, -6) node {$1$}; 
	\draw[blue] (-2.5, -6) node {$2$}; 
	\draw [purple,thick](0, -4) node (t)  {$\Rightarrow$};	
		
		\begin{scope}[shift={(4,0)}]
		
			{\tikzstyle{every node}=[draw ,circle,fill=white, minimum size=0.5cm,
			inner sep=0pt]
			\draw[blue,thick](-3,-2) node (s1)  {$s_1$};
			\draw[blue,thick](-3,-4) node (u)  {$u$};
			\draw[blue,thick](-3,-6) node (x)  {$x$};
			
			\draw[blue,thick](-1.5,-4) node (r)  {$r$};
			\draw[blue,thick](-1,-2.5) node (v1)  {$v_1$};
		}
		\path[draw,thick,black!60!green]
		(s1) edge node[name=la,pos=0.6, above] {\color{blue} $\tau\quad$} (u)
		(u) edge node[name=la,pos=0.6, above] {\color{blue}$1\quad$} (x)
		(r) edge node[name=la,pos=0.4, above] {\color{blue} $\tau\quad$} (v1);
		
		\draw[dashed, red, line width=0.5mm] (s1)--++(20:1cm); 
		\draw[dashed, red, line width=0.5mm] (s1)--++(160:1cm); 
		\draw[dashed, red, line width=0.5mm] (x)--++(340:1cm); 
		\draw[dashed, red, line width=0.5mm] (x)--++(200:1cm);
		\draw[dashed, red, line width=0.5mm] (r)--++(160:1cm);
		\draw[dashed, red, line width=0.5mm] (v1)--++(20:1cm); 
		
		\draw[blue] (-3.5, -2.1) node {$2$};
		\draw[blue] (-2.5, -2.1) node {$\D$};  
		\draw[blue] (-2.1, -4) node {$1$};  
		\draw[blue] (-0.6, -2.1) node {$1$}; 
		\draw[blue] (-3.5, -6) node {$\D$}; 
		\draw[blue] (-2.5, -6) node {$2$}; 
		\draw [purple,thick](0, -4) node (t)  {$\Rightarrow$};

		\end{scope}
		
			\begin{scope}[shift={(8,0)}]
		
		{\tikzstyle{every node}=[draw ,circle,fill=white, minimum size=0.5cm,
			inner sep=0pt]
			\draw[blue,thick](-3,-2) node (s1)  {$s_1$};
			\draw[blue,thick](-3,-4) node (u)  {$u$};
			\draw[blue,thick](-3,-6) node (x)  {$x$};
			
			\draw[blue,thick](-1.5,-4) node (r)  {$r$};
			\draw[blue,thick](-1,-2.5) node (v1)  {$v_1$};
		}
		\path[draw,thick,black!60!green]
		(s1) edge node[name=la,pos=0.6, above] {\color{blue} $\tau\quad$} (u)
		(u) edge node[name=la,pos=0.6, above] {\color{blue}$1\quad$} (x)
		(r) edge node[name=la,pos=0.4, above] {\color{blue} $1\quad$} (v1);
		
		\draw[dashed, red, line width=0.5mm] (s1)--++(20:1cm); 
		\draw[dashed, red, line width=0.5mm] (s1)--++(160:1cm); 
		\draw[dashed, red, line width=0.5mm] (x)--++(340:1cm); 
		\draw[dashed, red, line width=0.5mm] (x)--++(200:1cm);
		\draw[dashed, red, line width=0.5mm] (r)--++(160:1cm);
		\draw[dashed, red, line width=0.5mm] (v1)--++(20:1cm); 
		
		\draw[blue] (-3.5, -2.1) node {$2$};
		\draw[blue] (-2.5, -2.1) node {$\D$};  
		\draw[blue] (-2.1, -4) node {$\tau$};  
		\draw[blue] (-0.6, -2.1) node {$\tau$}; 
		\draw[blue] (-3.5, -6) node {$\D$}; 
		\draw[blue] (-2.5, -6) node {$2$}; 
		\draw [purple,thick](0, -4) node (t)  {$\Rightarrow$};

		\end{scope}
		
			\begin{scope}[shift={(12,0)}]
		
		{\tikzstyle{every node}=[draw ,circle,fill=white, minimum size=0.5cm,
			inner sep=0pt]
			\draw[blue,thick](-3,-2) node (s1)  {$s_1$};
			\draw[blue,thick](-3,-4) node (u)  {$u$};
			\draw[blue,thick](-3,-6) node (x)  {$x$};
			
			\draw[blue,thick](-1.5,-4) node (r)  {$r$};
			\draw[blue,thick](-1,-2.5) node (v1)  {$v_1$};
		}
		\path[draw,thick,black!60!green]
		(s1) edge node[name=la,pos=0.6, above] {\color{blue} $\tau\quad$} (u)
		(u) edge node[name=la,pos=0.6, above] {\color{blue}$1\quad$} (x)
		(r) edge node[name=la,pos=0.4, above] {\color{blue} $1\quad$} (v1);
		
		\draw[dashed, red, line width=0.5mm] (s1)--++(20:1cm); 
		\draw[dashed, red, line width=0.5mm] (s1)--++(160:1cm); 
		\draw[dashed, red, line width=0.5mm] (x)--++(340:1cm); 
		\draw[dashed, red, line width=0.5mm] (x)--++(200:1cm);
		\draw[dashed, red, line width=0.5mm] (r)--++(160:1cm);
		\draw[dashed, red, line width=0.5mm] (v1)--++(20:1cm); 
		
		\draw[blue] (-3.5, -2.1) node {$2$};
		\draw[blue] (-2.5, -2.1) node {$\D$};  
		\draw[blue] (-2.1, -4) node {$\tau$};  
		\draw[blue] (-0.6, -2.1) node {$2$}; 
		\draw[blue] (-3.5, -6) node {$\D$}; 
		\draw[blue] (-2.5, -6) node {$\tau$}; 
		\end{scope}
		
		\end{tikzpicture}
		-	  	\end{center}
	\caption{Three steps of Kempe changes}
	\label{f1}
\end{figure}
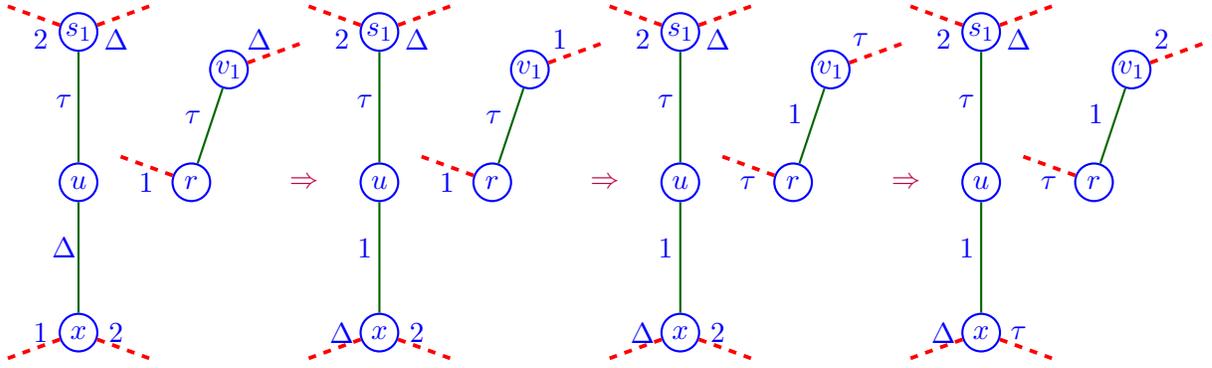

We first do three Kempe changes: Step 1: $(1,\Delta)$-swap(s) at both  $v_1$ and $x$ ($s_1$ and $r$ are $(1,\Delta)$-linked);   Step 2: a $(1,\tau)$-swap at $v_1$ (only changes the color on the edge $rv_1$); 
and Step 3:  $(2,\tau)$-swap(s) at both $x$ and $v_1$ ($s_1$ and $r$ are $(2,\tau)$-linked). See Figure~\ref{f1} for this sequence of changes. 

\medskip 
 
 

Note that 
Step 1 gives an $F$-stable coloring, Step 2 gives a $V(F-r)$-stable coloring, and Step 3 gives a stable coloring w.r.t. the new multifan obtained in Step 2.

%
%
%
%
%
%

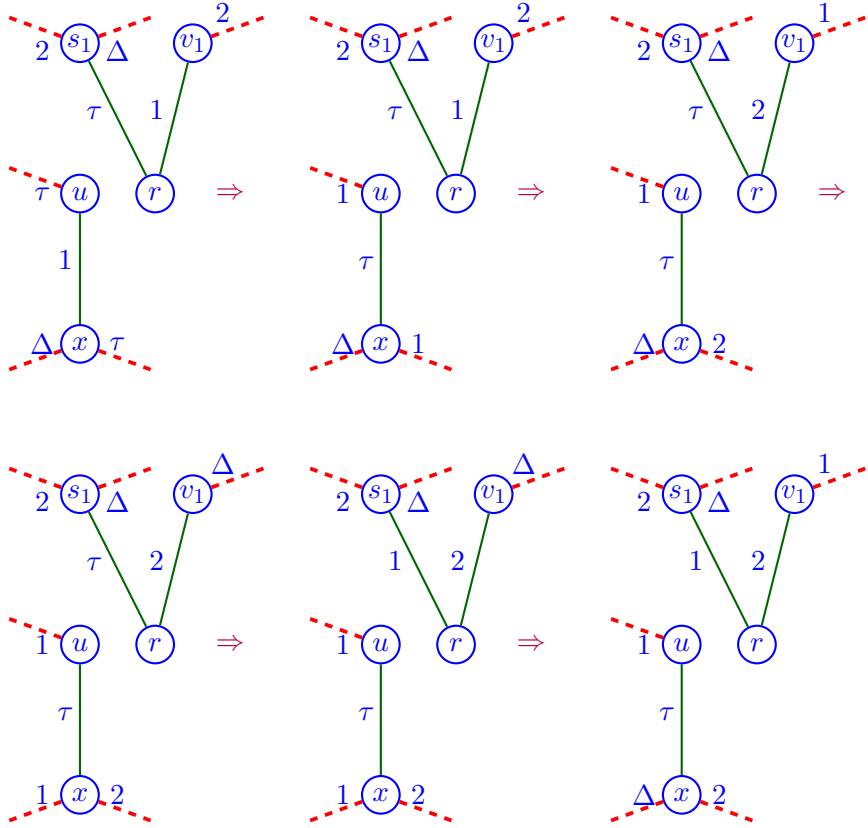
\begin{figure}[!htb]
 	\begin{center}
 		\begin{tikzpicture}[scale=1]
 		
 			{\tikzstyle{every node}=[draw ,circle,fill=white, minimum size=0.5cm,
 				inner sep=0pt]
 				\draw[blue,thick](-3,-2) node (s1)  {$s_1$};
 				\draw[blue,thick](-3,-4) node (u)  {$u$};
 				\draw[blue,thick](-3,-6) node (x)  {$x$};
 				
 				\draw[blue,thick](-2,-4) node (r)  {$r$};
 				\draw[blue,thick](-1.5,-2) node (v1)  {$v_1$};
 			}
 			\path[draw,thick,black!60!green]
 			(s1) edge node[name=la,pos=0.6, above] {\color{blue} $\tau\qquad$} (r)
 			(u) edge node[name=la,pos=0.6, above] {\color{blue}$1\quad$} (x)
 			(r) edge node[name=la,pos=0.4, above] {\color{blue} $1\quad$} (v1);
 			
 			\draw[dashed, red, line width=0.5mm] (s1)--++(20:1cm); 
 			\draw[dashed, red, line width=0.5mm] (s1)--++(160:1cm); 
 			\draw[dashed, red, line width=0.5mm] (x)--++(340:1cm); 
 			\draw[dashed, red, line width=0.5mm] (x)--++(200:1cm);
 			\draw[dashed, red, line width=0.5mm] (u)--++(160:1cm);
 			\draw[dashed, red, line width=0.5mm] (v1)--++(20:1cm); 
 			
 			\draw[blue] (-3.5, -2.1) node {$2$};
 			\draw[blue] (-2.5, -2.1) node {$\D$};  
 			\draw[blue] (-3.5, -4) node {$\tau$};  
 			\draw[blue] (-1.1, -1.6) node {$2$}; 
 			\draw[blue] (-3.5, -6) node {$\D$}; 
 			\draw[blue] (-2.5, -6) node {$\tau$}; 
 			
 			\draw [purple,thick](-1, -4) node (t)  {$\Rightarrow$};	
 		
 		\begin{scope}[shift={(4,0)}]
 			
 		{\tikzstyle{every node}=[draw ,circle,fill=white, minimum size=0.5cm,
 			inner sep=0pt]
 			\draw[blue,thick](-3,-2) node (s1)  {$s_1$};
 			\draw[blue,thick](-3,-4) node (u)  {$u$};
 			\draw[blue,thick](-3,-6) node (x)  {$x$};
 			
 			\draw[blue,thick](-2,-4) node (r)  {$r$};
 			\draw[blue,thick](-1.5,-2) node (v1)  {$v_1$};
 		}
 		\path[draw,thick,black!60!green]
 		(s1) edge node[name=la,pos=0.6, above] {\color{blue} $\tau\qquad$} (r)
 		(u) edge node[name=la,pos=0.6, above] {\color{blue}$\tau\quad$} (x)
 		(r) edge node[name=la,pos=0.4, above] {\color{blue} $1\quad$} (v1);
 		
 		\draw[dashed, red, line width=0.5mm] (s1)--++(20:1cm); 
 		\draw[dashed, red, line width=0.5mm] (s1)--++(160:1cm); 
 		\draw[dashed, red, line width=0.5mm] (x)--++(340:1cm); 
 		\draw[dashed, red, line width=0.5mm] (x)--++(200:1cm);
 		\draw[dashed, red, line width=0.5mm] (u)--++(160:1cm);
 		\draw[dashed, red, line width=0.5mm] (v1)--++(20:1cm); 
 		
 		\draw[blue] (-3.5, -2.1) node {$2$};
 		\draw[blue] (-2.5, -2.1) node {$\D$};  
 		\draw[blue] (-3.5, -4) node {$1$};  
 		\draw[blue] (-1.1, -1.6) node {$2$}; 
 		\draw[blue] (-3.5, -6) node {$\D$}; 
 		\draw[blue] (-2.5, -6) node {$1$}; 
 		
 		\draw [purple,thick](-1, -4) node (t)  {$\Rightarrow$};	
 		\end{scope}
 		
 		\begin{scope}[shift={(8,0)}]
 		
 		{\tikzstyle{every node}=[draw ,circle,fill=white, minimum size=0.5cm,
 			inner sep=0pt]
 			\draw[blue,thick](-3,-2) node (s1)  {$s_1$};
 			\draw[blue,thick](-3,-4) node (u)  {$u$};
 			\draw[blue,thick](-3,-6) node (x)  {$x$};
 			
 			\draw[blue,thick](-2,-4) node (r)  {$r$};
 			\draw[blue,thick](-1.5,-2) node (v1)  {$v_1$};
 		}
 		\path[draw,thick,black!60!green]
 		(s1) edge node[name=la,pos=0.6, above] {\color{blue} $\tau\qquad$} (r)
 		(u) edge node[name=la,pos=0.6, above] {\color{blue}$\tau\quad$} (x)
 		(r) edge node[name=la,pos=0.4, above] {\color{blue} $2\quad$} (v1);
 		
 		\draw[dashed, red, line width=0.5mm] (s1)--++(20:1cm); 
 		\draw[dashed, red, line width=0.5mm] (s1)--++(160:1cm); 
 		\draw[dashed, red, line width=0.5mm] (x)--++(340:1cm); 
 		\draw[dashed, red, line width=0.5mm] (x)--++(200:1cm);
 		\draw[dashed, red, line width=0.5mm] (u)--++(160:1cm);
 		\draw[dashed, red, line width=0.5mm] (v1)--++(20:1cm); 
 		
 		\draw[blue] (-3.5, -2.1) node {$2$};
 		\draw[blue] (-2.5, -2.1) node {$\D$};  
 		\draw[blue] (-3.5, -4) node {$1$};  
 		\draw[blue] (-1.1, -1.6) node {$1$}; 
 		\draw[blue] (-3.5, -6) node {$\D$}; 
 		\draw[blue] (-2.5, -6) node {$2$}; 
 		
 		\draw [purple,thick](-1, -4) node (t)  {$\Rightarrow$};	
 		\end{scope}
 		
 		\begin{scope}[shift={(0,-6)}]
 		
 		{\tikzstyle{every node}=[draw ,circle,fill=white, minimum size=0.5cm,
 			inner sep=0pt]
 			\draw[blue,thick](-3,-2) node (s1)  {$s_1$};
 			\draw[blue,thick](-3,-4) node (u)  {$u$};
 			\draw[blue,thick](-3,-6) node (x)  {$x$};
 			
 			\draw[blue,thick](-2,-4) node (r)  {$r$};
 			\draw[blue,thick](-1.5,-2) node (v1)  {$v_1$};
 		}
 		\path[draw,thick,black!60!green]
 		(s1) edge node[name=la,pos=0.6, above] {\color{blue} $\tau\qquad$} (r)
 		(u) edge node[name=la,pos=0.6, above] {\color{blue}$\tau\quad$} (x)
 		(r) edge node[name=la,pos=0.4, above] {\color{blue} $2\quad$} (v1);
 		
 		\draw[dashed, red, line width=0.5mm] (s1)--++(20:1cm); 
 		\draw[dashed, red, line width=0.5mm] (s1)--++(160:1cm); 
 		\draw[dashed, red, line width=0.5mm] (x)--++(340:1cm); 
 		\draw[dashed, red, line width=0.5mm] (x)--++(200:1cm);
 		\draw[dashed, red, line width=0.5mm] (u)--++(160:1cm);
 		\draw[dashed, red, line width=0.5mm] (v1)--++(20:1cm); 
 		
 		\draw[blue] (-3.5, -2.1) node {$2$};
 		\draw[blue] (-2.5, -2.1) node {$\D$};  
 		\draw[blue] (-3.5, -4) node {$1$};  
 		\draw[blue] (-1.1, -1.6) node {$\Delta$}; 
 		\draw[blue] (-3.5, -6) node {$1$}; 
 		\draw[blue] (-2.5, -6) node {$2$}; 
 		
 		\draw [purple,thick](-1, -4) node (t)  {$\Rightarrow$};	
 		\end{scope}

 			\begin{scope}[shift={(4,-6)}]
 		
 		{\tikzstyle{every node}=[draw ,circle,fill=white, minimum size=0.5cm,
 			inner sep=0pt]
 			\draw[blue,thick](-3,-2) node (s1)  {$s_1$};
 			\draw[blue,thick](-3,-4) node (u)  {$u$};
 			\draw[blue,thick](-3,-6) node (x)  {$x$};
 			
 			\draw[blue,thick](-2,-4) node (r)  {$r$};
 			\draw[blue,thick](-1.5,-2) node (v1)  {$v_1$};
 		}
 		\path[draw,thick,black!60!green]
 		(s1) edge node[name=la,pos=0.6, above] {\color{blue} $1\qquad$} (r)
 		(u) edge node[name=la,pos=0.6, above] {\color{blue}$\tau\quad$} (x)
 		(r) edge node[name=la,pos=0.4, above] {\color{blue} $2\quad$} (v1);
 		
 		\draw[dashed, red, line width=0.5mm] (s1)--++(20:1cm); 
 		\draw[dashed, red, line width=0.5mm] (s1)--++(160:1cm); 
 		\draw[dashed, red, line width=0.5mm] (x)--++(340:1cm); 
 		\draw[dashed, red, line width=0.5mm] (x)--++(200:1cm);
 		\draw[dashed, red, line width=0.5mm] (u)--++(160:1cm);
 		\draw[dashed, red, line width=0.5mm] (v1)--++(20:1cm); 
 		
 		\draw[blue] (-3.5, -2.1) node {$2$};
 		\draw[blue] (-2.5, -2.1) node {$\D$};  
 		\draw[blue] (-3.5, -4) node {$1$};  
 		\draw[blue] (-1.1, -1.6) node {$\Delta$}; 
 		\draw[blue] (-3.5, -6) node {$1$}; 
 		\draw[blue] (-2.5, -6) node {$2$}; 
 		
 		\draw [purple,thick](-1, -4) node (t)  {$\Rightarrow$};	
 		\end{scope}

 		\begin{scope}[shift={(8,-6)}]
 		
 		{\tikzstyle{every node}=[draw ,circle,fill=white, minimum size=0.5cm,
 			inner sep=0pt]
 			\draw[blue,thick](-3,-2) node (s1)  {$s_1$};
 			\draw[blue,thick](-3,-4) node (u)  {$u$};
 			\draw[blue,thick](-3,-6) node (x)  {$x$};
 			
 			\draw[blue,thick](-2,-4) node (r)  {$r$};
 			\draw[blue,thick](-1.5,-2) node (v1)  {$v_1$};
 		}
 		\path[draw,thick,black!60!green]
 		(s1) edge node[name=la,pos=0.6, above] {\color{blue} $1\qquad$} (r)
 		(u) edge node[name=la,pos=0.6, above] {\color{blue}$\tau\quad$} (x)
 		(r) edge node[name=la,pos=0.4, above] {\color{blue} $2\quad$} (v1);
 		
 		\draw[dashed, red, line width=0.5mm] (s1)--++(20:1cm); 
 		\draw[dashed, red, line width=0.5mm] (s1)--++(160:1cm); 
 		\draw[dashed, red, line width=0.5mm] (x)--++(340:1cm); 
 		\draw[dashed, red, line width=0.5mm] (x)--++(200:1cm);
 		\draw[dashed, red, line width=0.5mm] (u)--++(160:1cm);
 		\draw[dashed, red, line width=0.5mm] (v1)--++(20:1cm); 
 		
 		\draw[blue] (-3.5, -2.1) node {$2$};
 		\draw[blue] (-2.5, -2.1) node {$\D$};  
 		\draw[blue] (-3.5, -4) node {$1$};  
 		\draw[blue] (-1.1, -1.6) node {$1$}; 
 		\draw[blue] (-3.5, -6) node {$\D$}; 
 		\draw[blue] (-2.5, -6) node {$2$}; 
 	 		 		\end{scope}

 		\end{tikzpicture}
 		-	  	\end{center}
 	\caption{Five steps of Kempe changes}
 	\label{f2}
 \end{figure}

We then  color $rs_1$ by $\tau$ and uncolor $s_1u$ to give a coloring $\phiv'$, which is followed by  
5 Kempe changes:  $ux$:  $1 \rightarrow \tau$,   $(1,2)$-swap(s) at both $x$
and $v_1$ ($s_1$ and $u$ are $(1,2)$-linked),   $(1,\Delta)$-swap(s) at both $x$ and $v_1$ ($s_1$ and $u$ are $(1,\D)$-linked), 
 a $(1,\tau)$-swap on the $(1,\tau)$-chain containing $s_1r$, and   $(1,\Delta)$-swap(s) at both $x$ and $v_t$ ($s_1$ and $u$ are $(1,\D)$-linked). 
 Since every recoloring is a Kempe change, the final coloring is in $\CC^{\D}(G-s_1u)$. 
 See Figure~\ref{f2} for this sequence of changes.  
 
 \medskip 

Under the current coloring, we have $P_{s_1}(1,2)=s_1rv_1$. On the other hand, since $s_1u$ is uncolored and $1$ and $2$ are missing at $u$ and $s_1$ respectively,  $P_{s_1}(1, 2) = P_u(1, 2)$, giving a contradiction.

{\bf \noindent Case 2. $t\ge 2$}. 

\smallskip 

Let $\varphi(rv_t)=\pbar(v_{t-1})=\tau_t$. We may assume 
that $v_{t-1}$ and $r$ are $(1, \tau_t)$-linked. Otherwise, the 
$(1, \tau_t)$-swap at $v_{t-1}$  gives an $F$-stable coloring that contradicts 
Claim~\ref{cla:noB1}. We do a $(1,\tau_t)$-swap at $x$ and get an $F$-stable coloring. 
By Lemma~\ref{lem:rs1}, 
$r\in P_{s_1}(\tau_t,\Delta)=P_{v_t}(\tau_t,\Delta)$. We then do 
 $(\tau_t,\Delta)$-swaps at both $x$ and $v_{t-1}$ and get an $F$-stable coloring $\phiv'$. Note that
 $\phiv'(ux) = \tau_t$ and $\D \in \pbar'(x) \cap \pbar'(v_{t-1})\cap\pbar'(v_t)$. 

By Lemma~\ref{lem:rs1}, $r\in P_{s_1}(\tau, \Delta)$. 
We  claim that $P_{s_1}(\tau, \Delta)=P_x(\tau, \Delta)$. 
Suppose to the contrary that these two paths are disjoint. 
 If $P_x(\tau,\Delta)$ is also disjoint from  $P_{v_t}(\tau,\Delta)$, we 
do the following sequence of five Kempe changes:  the $(\tau, \Delta)$-swap at $x$,   the  $(1,\Delta)$-swap at $v_t$ ($s_1$ and $r$ are $(1,\Delta)$-linked),
the $(1,\tau_t)$-swap on the $(1,\tau)$-chain containing $ux$, the $(2,\Delta)$-swap at $x$, and the $(1,\Delta)$-swap at $x$. 
Except the Kempe change  that the $(2,\Delta)$-swap at $x$ may possibly change the colors on two edges of $F$,
all other changes are $F$-stable. Thus the final resulting coloring is $V(F)$-stable. 
 Under the current coloring, $P_{s_1}(\tau, \D) = s_1ux$ that does not contain vertex $r$, giving a contradiction to Lemma~\ref{lem:rs1}. 
%
%

%

Under the assumption that $P_{s_1}(\tau, \Delta) \ne P_x(\tau, \Delta)$, by the argument above, we assume then that 
 $P_x(\tau,\Delta) = P_{v_t}(\tau,\Delta)$. We 
do the $(\tau, \Delta)$-swap at $x$ that is also the $(\tau, \D)$-swap at $v_t$ to get an $F$-stable coloring.  Note that 
$\D$ is no longer missing at $x$ unless $\tau$ is also previously missing at $x$. 
Since $s_1$ and $r$ are $(1,\Delta)$-linked,  we do a $(1,\Delta)$-swap at $v_{t-1}$ to get an $F$-stable coloring,  
and do a shifting from $v_1$ to $v_{t-1}$, which give a $V(F-r)$-stable coloring. Denote the corresponding new multifan by $F^*$. Since $s_1$ and $r$ are $(\tau,\Delta)$-linked, 
we do a $(\tau, \Delta)$-swap at both $x$ and $v_t$ to get an $F^*$-stable coloring such that $\D$ is missing at $x$.  
Since $r\in P_{s_1}(\tau_t,\Delta)=P_{v_t}(\tau_t,\Delta)$ by Lemma~\ref{lem:rs1}, 
we  do the $(\Delta, \tau_t)$-swap at $x$, which does not affect the multifan.  
Since $\phiv^*(s_1u) =\tau =\pbar^*(r)$ and $\phiv^*(ux) =\D \in \pbar^*(s_1)$, $(r, rs_1, s_1, s_1u, u, ux, x)$ is 
a Kierstead path. But, $2\in \pbar^*(s_1)\cap \pbar^*(x)$, giving a contradiction. 
Therefore  $r\in P_{s_1}(\tau, \Delta, \phiv')=P_x(\tau, \Delta, \phiv')$.
 
We do a $(\tau, \Delta)$-swap at both $v_{t-1}$
and $v_t$. Under the new coloring, $\tau$ is missing at both $v_{t-1}$ and $v_t$. We may assume that  $v_t$ and $r$  are $(1,\tau)$-unlinked by doing the $A$-shifting from $v_1$
to $v_{t-1}$ if necessary. 
Thus we do a $(1,\tau)$-swap at $v_t$. 
Denote the new coloring by $\varphi^*$. If $\varphi^*(s_1u)=\tau$, we do a $(1,\tau_t)$-swap on the $(1,\tau_t)$-chain containing $ux$ and then do a $(1,\Delta)$-swap at $x$. This gives an A-type $\tau$-sequence. 
Thus, $\varphi^*(s_1u)=1$.  We do a $(1,\tau_t)$-swap on the $(1,\tau_t)$-chain containing $ux$
and then do a $(1,\Delta)$-swap at both $x$ and $v_t$. This gives back to Case 1 with $v_t$ in the place of 
$v_1$ and $\tau_t$ in the place of $\tau$. 
\qed

\subsection{Proof of Lemma~\ref{lem:good-coloring}} 

Assume to  the contrary that  there exists a vertex $u\in N(s_1)\cap N(x)$ with $u\not\in N(r)\setminus N_\Delta(r)$ such that    
$2,\Delta \in \pbar(x)$ and $\varphi(s_1u)=1$. Note that $u\ne r$, as every neighbor of $r$ has degree at least $\Delta-1$ in $G$ while $d(x) \le \Delta-3$. Thus $u\not\in N[r]\setminus N_\Delta(r)$. 
Let $\varphi'(ux)=\tau$. Clearly, $\tau \ne 1$. 

Since $F$ is a maximum multifan,  $r\in P_{s_1}( \tau,\D)$ and $r\in P_{s_1}(2, \tau)$ by Lemma~\ref{lem:rs1}.  We claim that
$P_{s_1}(\tau,\D) = P_x(\tau,\D)$ and $P_{s_1}(2, \tau) = P_x(2, \tau)$. Otherwise, say $P_x(\tau,\D)$ and $P_{s_1}(\tau,\D)$ are disjoint. We do a
$(\tau,\D)$-swap at $x$ and get an $F$-stable coloring $\phiv'$.  Since 
$\phiv'(s_1u) =\phiv(s_1u) =1 \in \pbar'(r)$ and $\phiv'(ux) =\D \in \pbar'(s_1)$, $K=(r, rs_1,s_1,s_1u, u, ux,x)$ is a Kierstead path and 
$2\in \pbar'(s_1)\cap \pbar'(x)$,  contradicting $V(K)$ being elementary  (Lemma~\ref{Lemma:kierstead path1}).

 We claim that  $\tau \notin \pbar(F)$. Otherwise, $P_{s_1}(\tau, \D) = P_{\pbar^{-1}_F(\tau)}(\tau, \D)$ if $\tau$ is $2$-inducing and 
 $P_{s_1}(2, \tau) = P_{\pbar^{-1}_F(\tau)}(2, \tau)$ if $\tau$ is $\D$-inducing. In either case, we get a contradition to the previous claim. 
Since the multifan $F$ is maximum,   there is a unique $\tau$-sequence  $(v_1, \ldots, v_t)$ by Lemma~\ref{lem:tau-seq}. 
Since $r\in P_{s_1}(2,\tau)=P_x(2,\tau)$ and $r\in P_{s_1}(\tau,\Delta)=P_x(\tau, \Delta)$, $rv_1 \in P_x(2, \tau)$ and $rv_1\in P_x(\tau, \D)$. 

If the $\tau$-sequence is of type A,  we do the A-shifting and get an $F$-stable coloring $\varphi'$, and under this coloring $P_x(\D, \tau, \phiv') \ne P_{s_1}(\D, \tau, \phiv')$. But, $\pbar'(x) =\pbar(x) \supseteq \{2, \D\}$ and $\phiv'(s_1u) = \phiv(s_1u) =1$, giving a contradiction.

Suppose then that  the $\tau$-sequence is of type B:  $\pbar(v_t)=\gamma\in \pbar(F)$.  
If $\gamma=1$, we do the B-shifting and get a $V(F-r)$-stable coloring. Note that $F$ is also a multifan w.r.t. the new coloring and that $\tau$ and $\D$ are missing at $r$ and $s_1$, respectively.  We then do a $(\tau, \D)$-swap at $x$ and get an $F$ stable coloring $\phiv'$ w.r.t. the previous coloring. Note that $F$ is a multifan w.r.t. $rs_1$ and $\phiv'$, $\pbar'(r) = \tau$, $\{2, \tau\} \subseteq \pbar'(x)$ and $\phiv'(ux) =\D$,    showing a contradiction to Lemma~\ref{lem:good-coloring0}. 
So, $\gamma \ne 1$, say $\gamma$ is 2-inducing. 
Let $\varphi'=\varphi/ P_{v_t}(1,\gamma,\varphi)$. 
If $s_1u\not \in P_{v_t}(1,\gamma,\varphi)$, the argument turns back to  $\pbar(v_t)=\gamma=1$ case, which we just settled.  
Thus, we assume that $\varphi'(s_1u)=\gamma$. 
We do a shifting from $v_1$ to $v_t$. Then as $s_1$ and $r$ are $(\tau,\Delta)$-linked, 
we do a $(\tau,\Delta)$-swap at $x$. Now up to exchanging the role of 1 and $\tau$, we have an $V(F)$-stable coloring 
$\varphi''$ such that  $\varphi''(s_1u)=\gamma$, $\varphi''(ux)=\Delta$
and $1,2\in \pbar''(x)$, and $\gamma$ is a 2-inducing color of $F$ w.r.t. $\varphi''$. 
Let $\varphi^*=\varphi''/P_x(1,\gamma)$. If $s_1u\not\in P_x(1,\gamma,\varphi'')$, then $P_{s_1}(\gamma, \Delta)=s_1ux$, 
showing a contradiction to Lemma~\ref{thm:vizing-fan2}~\eqref{thm:vizing-fan2-a}. 
If  $s_1u\in P_x(1,\gamma,\varphi'')$, then $K=(r, rs_1,s_1,s_1u, u, ux,x)$ is a Kierstead path, but 2 is missing at both $s_1$
and $x$, showing a contradiction to Lemma~\ref{Lemma:kierstead path1}.

Thus the $\tau$-sequence is of type C:  $\pbar(v_t)= \pbar(v_{i-1})=\tau_i$  for some $i\in [2,t-1]$ and some $\tau_i\in [1,\D]\setminus \pbar(F)$.  We first do a $(1,\Delta)$-swap at $x$.
One of $v_{i-1}$
and $v_t$ is   $(\tau_i,\Delta)$-unlinked with $s_1$.  
We may assume that 
$v_t$ and $s_1$ are $(\tau_i,\Delta)$-unlinked (the proof for the other case is similar).
By Lemma~\ref{lem:rs1}, we have 
$r\not\in P_{v_t}(\tau_i,\Delta)$. 
We first do a $(\tau_i, \Delta)$-swap at $v_t$ and then a $(1,\Delta)$-swap at both $x$
and $v_t$. This converts the problem back to the   type B $\tau$-sequence case.  
\qed

\subsection{Proof of Lemma~\ref{lem:good-coloring1}}

For (i), let $1, \tau, \lambda\in \pbar(x)$ be three distinct colors.
Suppose to the contrary that there does not exist a $V(F)$-stable coloring such that 
both $2, \D$ are missing at $x$.   


 We claim  $\{\tau, \lambda\}\cap \pbar(F) = \emptyset$. Otherwise, say  $\tau \in \pbar(F)$.  We do a $(1,2)$-swap at $x$, and then do $(\tau,1)-(1,\Delta)$-swaps at $x$
to get a new coloring $\phiv'$.  Clearly, $2, \D\in \pbar'(x)$.   Since $1\in \pbar(r)$ and $ \tau, \D\in \pbar(F)$,  $\phiv'$ is $F$-stable, giving a contradiction. 

 We note that none of the following four Kempe changes   $\varphi/P_x(2,\tau)$,  $\varphi/P_x(\tau,\Delta)$,  $\varphi/P_x(2,\lambda)$, and $\varphi/P_x(\lambda, \D)$ is  $V(F)$-stable. As otherwise, say that $\varphi/P_x(2,\tau)$ is $V(F)$-stable, 
 we do a $(2,\tau)$-swap at $x$ and then a $(1,\Delta)$-swap at $x$ in getting a desired coloring.  Applying  Lemma~\ref{lem:rs1},  we get 
 $r\in P_{s_1}(2,\tau)=P_x(2,\tau)$, $r\in P_{s_1}(\tau,\Delta)=P_x(\tau,\Delta)$, $r\in P_{s_1}(2,\lambda)=P_x(2,\lambda)$, and $r\in P_{s_1}(\lambda,\Delta)=P_x(\lambda,\Delta)$.  
 By Lemma~\ref{lem:tau-seq2} (v),  the $\tau$-sequence $(v_1, \ldots, v_t)$
 is of type $B$ with $\pbar(v_t)=1$ and   the $\lambda$-sequence $(w_1, \ldots, w_k)$
 is of type $B$ with $\pbar(w_k)=1$. 
 
 If these two sequences are disjoint, then $x$
 is $(1,2)$-linked with at most one of $v_t$ and $w_k$.  Assume, without loss of generality, that $x$
 and $v_t$ are $(1,2)$-unlinked. 
 We do a $(1,2)$-swap at $x$ to get an $F$-stable coloring, and then do a B-shifting from $v_1$ to $v_t$ to get a $V(F-r)$-stable coloring. (Note that $F$ is still a maximum multifan w.r.t. the new coloring.) 
 Next, we do a $(\tau,\Delta)$-swap at $x$ and get an $F$-stable coloring $\phiv'$ w.r.t. the previous coloring. So, $\phiv'$ is a $V(F-r)$-stable coloring such that $2, \D\in \pbar'(x)$. By switching colors $\tau$ and $1$ for the entire graph, we get a $V(F)$-stable coloring such that $2, \D$ are missing at $x$, a contradiction.  Therefore, the $\tau$-sequence and the $\lambda$-sequence have some overlap.  Assume that $v_i =w_j$ is the first common vertex of the two sequence.  Then, the two sequences are identical after  this vertex. 
 
 If both $i, j$ are at least two, then $\pbar(v_{i-1}) = \pbar(w_{j-1})$, name it $\gamma$. By the definition of 
 $\tau$-sequence, $\gamma\in [1,\Delta]\setminus \pbar(F)$.  
 Since $F$ is maximum, $r\in P_{s_1}(\gamma, \D)$ by Lemma~\ref{lem:rs1}. One of  $v_{i-1}$ and $w_{j-1}$, say $v_{i-1}$, is not  on $P_{s_1}(\gamma, \D)$. We do a $(\gamma, \D)$-swap at $v_{i-1}$ and get an $F$-stable coloring $\phiv'$.  But, the $\tau$-sequence ends with a vertex missing color $\D$ rather than $1$, giving a contradiction to Lemma~\ref{lem:tau-seq2} (v).

Assume without loss of generality $j=1$, i.e.,  $\lambda=\pbar(v_{i-1})$ and $w_1, \ldots, w_k$ is the same as $v_i, \ldots, v_t$. 
 We first do a $(1,2)$-swap at both $x$ and $v_t$.
 One of $x$ and $v_{i-1}$ is $(1,\lambda)$-unlinked with $r$. 
 If $x$ and $r$ are $(1,\lambda)$-unlinked, we do $(\lambda,1)-(1,\Delta)$-swaps at $x$ to get a desired coloring. 
 If $v_{i-1}$ and $r$ are $(1,\lambda)$-unlinked, we do a $(1,\lambda)$-swap at $v_{i-1}$
 and then do a shifting from $v_1$ to $v_{i-1}$. Next, we do a $(\tau,\Delta)$-swap at $x$ and get a $V(F-r)$-coloring $\phiv'$.  Switching colors $1$ and $\tau$ for the entire graph, we get a $V(F)$-stable coloring such that $2, \D$ are missing at $x$, a contradiction.


 For (ii)(a),  assume that $2, \D\in \pbar(x)$ and $\varphi(ux)=\tau$.
 We first note that  $\tau \notin \pbar(F)$. Otherwise,  since $\tau$ is in $\pbar(r)$,   $2$-inducing, or $\D$-inducing, we assume that $\tau$ is not $\D$-inducing. 
 We do  $(\Delta, \tau)-(\tau,1)$-swaps at $x$ and get a coloring $\phiv'$.  
 Since $x$  and $s_1$ are  $(\tau, \D)$-unlinked and $r$ and  $\pbar_F^{-1}(\tau)$ are $(1,\tau)$-linked,   $\phiv'$ is  $V(F)$-stable. But  
 $\varphi'(ux)=\Delta$ and $1,2\in \pbar(x)$, showing a contradiction to Lemma~\ref{lem:good-coloring0}.

 By Lemma~\ref{lem:tau-seq}, there is a unique $\tau$-sequence $(v_1, \dots, v_t)$. Assume without loss of generality that $\phiv(s_1u)$ is not $\D$-inducing. We will show that there exists a $V(F)$-stable coloring $\phiv'$ satisfying the following conditions: 
 \begin{equation} \label{eq:lem4.3-ii-a}
 1\in\pbar'(x), \, \phiv'(ux) =\D, \mbox{ and $\phiv'(su)$ is not a $\D$-inducing color}.
 \end{equation}  
 Suppose to the contrary that 
 such a $V(F)$-stable coloring does not exist. 
  
 Since $F$ is a maximum multifan, $r\in P_{s_1}(\tau, \D)$. We claim $x$ and $s_1$ are $(\tau, \D)$-linked, and so 
 $rv_1\in P_{s_1}(\tau, \D, \phiv) = P_x(\tau, \D, \phiv)$.   Otherwise, we do the  $(\tau, \D)$-swap and the $(1, 2)$-swap at $x$ and get a coloring $\phiv'$.  Clearly,  $\phiv'$ is an $F$-stable coloring and satisfies (\ref{eq:lem4.3-ii-a}),  a contradiction.   By considering the type of the $\tau$-sequence, we will find a $V(F)$-stable coloring such that $x$ and $s_1$ are not $(\tau, \D)$-linked, which in turn gives a contradiction.

 {\flushleft  Suppose that the $\tau$-sequence is of type A}.   We do the $A$-shifting on the $\tau$-sequence to obtain an $F$-stable coloring $\phiv'$.  
 Clearly $\phiv'$ is $F$-stable  and  either $P_{s_1}(\tau, \D,\phiv') =P_{v_1}(\tau, \D, \phiv')$ or 
 $P_{x}(\tau, \D,\phiv') =P_{v_1}(\tau, \D, \phiv')$, but not both.   So, $x$ and $s_1$ are not $(\tau, \D)$-linked, a contradiction.  
 
 {\flushleft Suppose that  the $\tau$-sequence is of type B}. Let $\pbar(v_t)=\gamma$ that is a color in $\pbar(F)$. If $\gamma \ne \D$, we do the
  $(1,\gamma)$-swap at $v_t$ and then do a B-shifting from $v_1$ to $v_t$ (when $\gamma=2$, the color 2 missing at x might be changed to 1 after this swap). The resulting coloring $\phiv'$ is $V(F -r)$-stable (with $1$ being replaced by $\tau$ at $r$) and  $P_{s_1}(\tau,\Delta, \phiv')\ne P_{x}(\tau,\Delta, \phiv')$. Note that we either have $\phiv'(s_1u) = \phiv(s_1u)$ or $\phiv'(s_1u) =1$, so in either case $\phiv'(s_1u)$ is still not a $\D$-inducing color.  As $s_1$ and $r$ are $(\tau,\Delta)$-linked, we do a $(\tau,\D)$-swap at $x$. 
  By switching colors $1$ and $\tau$ for the entire graph, we get a $V(F)$-stable coloring satisfying (\ref{eq:lem4.3-ii-a}), a contradiction. 
  Thus $\gamma = \D$. We do $(1, \D)$-swaps at both $x$ and $v_t$.  (It may be only one swap if $x$ and $v_t$ are $(1, \D)$-linked.)  Note that $1$ is missing at both $x$ and $v_t$ now.  Do the $(1, 2)$-swap at $v_t$ and the $(1, \D)$-swap at $x$. Clearly, the resulting coloring $\phiv'$  is $F$-stable with  $\phiv'(s_1u) =\phiv(s_1u)$, and $ \pbar'(v_t)=2$. This turns the problem back to the previous case where $\gamma\ne \Delta$. 
  
  {\flushleft Suppose that  the $\tau$-sequence is of type C}. That is, $\pbar(v_t) = \pbar(v_{i-1})=\tau_i$ for some $i\in [2,t-1]$ and some $\tau_i\in [1,\Delta]\setminus \pbar(F)$. Since one of $v_{i-1}$ and $v_t$, say $v_t$,  is $(1, \tau_i)$-unlinked with $r$, we do $(1, \tau_i)$-swap at $v_t$ and get an $F$-stable coloring $\phiv'$. Note that $\phiv'(s_1u)$ is either $\phiv(s_1u)$ or $1$, and so it is still not $\D$-inducing. Under coloring
  $\phiv'$, the $\tau$-sequence $(v_1, \dots, v_t)$ is of type B with $\pbar'(v_t) = 1$, which was previously settled.


 For (ii)(b), by  (ii)(a), we assume, without loss of generality,  that there is a $V(F)$-stable coloring  $\varphi$ such that $\varphi(ux)=\Delta$, $1\in \pbar(x)$, and $\varphi(s_1u)$ is not 
 $\Delta$-inducing.   We show that there is an $V(F)$-stable coloring $\varphi'$
 such that $\varphi'(s_1u)=1$ and $\pbar'(x) \cap \{2, \D\} \ne \emptyset$.  
 Let $\varphi(s_1u)=\tau$. 
 
 Assume first that $\tau \in \pbar(F)$.  If  $\tau = 1$, then $K=(r, rs_1,s_1, s_1u, u, ux, x)$
 is a Kierstead path since $\phiv(s_1u) \in \pbar(r)$ and $\phiv(ux) =\D \in \pbar(s_1)$. But $1\in \pbar(s_1)\cap \pbar(x)$ gives a contradiction to Lemma~\ref{Lemma:kierstead path1}. Thus $\tau\ne 1$, and so is $2$-inducing
 as it is not $\D$-inducing by (ii) (a). 
We do a $(1,\tau)$-swap at $x$ and get 
 an  $F$-stable coloring $\varphi^*$.  
 We have either $\phiv^*(s_1u) = \phiv(s_1u) =\tau$ or $\phiv^*(s_1u) = 1$. 
 If $\varphi^*(s_1u)=\tau$, then $P_{s_1}(\tau,\Delta)=s_1ux$, contradicting Lemma~\ref{thm:vizing-fan2}~\eqref{thm:vizing-fan2-a}.
 Thus, $\varphi^*(s_1u)=1$.
 Since $\varphi(ux)=\Delta$ and $\tau$ is 
  $2$-inducing, 
  we do a $(\tau,\Delta)$-swap at $x$ 
 and  get a desired $V(F)$-stable coloring.

 We now assume  $\tau\in [1,\Delta]\setminus \pbar(F)$. We first do a $(1,2)$-swap at $x$ and get an $F$-stable coloring
 $\phiv'$.  
Notice that  $\varphi'(ux)=\Delta$ and $2\in \pbar'(x)$. 
Let $v_1, \ldots, v_t$ be the $\tau$-sequence by Lemma~\ref{lem:tau-seq}. 
We claim that for any   $V(F)$-stable and $\{2, \tau, \D\}$-avoiding coloring $\phiv''$ w.r.t. $\phiv'$, 
$P_{x}(2,\tau) =P_{s_1}(2, \tau)$. Otherwise, we do a $(2, \tau)$-swap at $x$ and get an $F$-stable coloring $\varphi''$. 
Under the new coloring, since $\phiv''(ux) = \phiv'(ux) =\D$ and $\phiv''(s_1u) = \phiv'(s_1u) = \tau$,
we have $P_{s_1}(\tau, \D) = s_1 ux$, contradicting $r\in P_{s_1}(\tau, \D)$ (Lemma~\ref{lem:rs1}). 
Applying Lemma~\ref{lem:tau-seq2} (iv) with $2$ replacing $\D$,  we see that  the $\tau$-sequence is of type B and 
 $\pbar'(v_t)\in \{2,\Delta\}$. 
 
 If $\pbar'(v_t)=2$, since $r\in P_{s_1}(2, \tau, \phiv')=P_{x}(2,\tau, \phiv')$, 
 we do the $(2,\tau)$-swap at $v_t$ and get an $F$-stable  $\phiv''$. 
 Note that  $\phiv''(ux) =\phiv'(ux) = \D$ and $\phiv''(s_1u) = \phiv'(s_1u) =\tau$ since $s_1u\in P_{s_1}(2, \tau)$. For the same reason
 above, we still have $r\in P_{s_1}(2, \tau, \phiv'') = P_x(2, \tau, \phiv'')$. 
 Under coloring $\phiv''$, the $\tau$-sequence $(v_1, \dots, v_t)$ is of type A. 
 We do the $A$-shifting from $v_1$ to $v_t$ and get an $F$-stable coloring $\phiv'''$. 
 Under the new coloring $\phiv'''$,  we still have 
 $\phiv'''(ux) =\D$ and $\phiv'''(s_1u) =\tau$,  but $ P_{s_1}(2, \tau, \phiv''') \ne P_{x}(2,\tau, \phiv''')$ gives a contradiction. 
 
  If $\pbar'(v_t)=\Delta$,
 then we do the  $(1,\Delta)$-swap at $v_t$ and do the  $B$-shifting from $v_1$
 to $v_t$. Denote the new coloring by $\varphi''$. Notice that $\phiv''(f) =\phiv'(f)$ for every $f\in E(F)$, $\pbar''(v) =\pbar''(v)$ for every $v\in V(F-r)$, and $\tau\in \pbar''(r)$. 
 If  $ux \not\in P_{v_t}(1,\Delta,\varphi')$, then $\varphi''(ux)=\Delta\in\pbar''(s_1)$ and $\phiv''(s_1u) =\phiv'(s_1u) =\tau \in \pbar''(r)$. Thus  $K=(r,rs_1,s_1,s_1u,ux,x)$ is a Kierstead path but 
 $2\in \pbar''(s_1)\cap \pbar''(x)$, showing a contradiction to Lemma~\ref{Lemma:kierstead path1}. 
 Thus, $ux \in P_{v_t}(1,\Delta,\varphi')$, which in turn gives $\varphi''(ux)=1$. Then, by switching $1$ and $\tau$ for the entire graph, we get a desired $V(F)$-stable coloring.

 For the moreover part of (ii)(b), we assume $\varphi(s_1x)=1$ and $2\in \pbar(x)$ and 
 show that there is a $V(F)$-stable coloring $\phiv^*$
 such that $\varphi^*(s_1x)=1$ and $2, \Delta\in \pbar^*(x)$.
Let $\tau\in \pbar(x)\setminus\{2\}$. 

Suppose $\tau \in \pbar(F)$. If $\tau$ is not $\D$-inducing, we simply do a $(\tau,\Delta)$-swap at $x$ to get a desired $V(F)$-stable coloring.  If $\tau$ is $\Delta$-inducing, we do  $(2,1)-(1,\Delta)$-swaps at $x$ to get an $F$-stable coloring such that $\tau$ is still
$\D$-inducing,  and then do a $(2,\tau)$-swap at $x$ to get a desired $V(F)$-stable coloring. 

Suppose $\tau\in [1,\Delta]\setminus \pbar(F)$. We first do a $(1,2)$-swap
at $x$ and get an $F$-stable coloring $\phiv'$.  We now have 
$1,\tau \in \pbar'(x)$. Since $s_1u \in P_{s_1}(1, 2)=P_r(1,2)$, we have $\phiv'(s_1u) = \phiv(s_1u) =1$. 

Since $F$ is a maximum multifan, by Lemma~\ref{lem:rs1}
we have $r\in P_{s_1}(2, \tau)$ and $r\in P_{s_1}(\tau, \D)$. 
We claim that for any $V(F)$-stable and $(1, \tau)$-avoiding coloring $\phiv''$, we have
$P_x(\tau,\Delta, \phiv'') = P_{s_1}(\tau, \D, \phiv'')$.  Otherwise, 
a $(\tau, \D)$-swap at $x$ and then a $(1, 2)$-swap at $x$ give a $V(F)$-stable coloring $\phiv''$ such that
$2, \D \in \pbar''(x)$ and $\phiv''(s_1u) = \phiv'(s_1u) =1$ since $s_1u\in P_{s_1}(1, 2) = P_r(1, 2)$. 
 By Lemma~\ref{lem:tau-seq2} (v), 
the $\tau$-sequence 
 $(v_1, \ldots, v_t)$ is of type B
such that  $\pbar'(v_t)=1$.

Since $d(x)\le \Delta-3$,  there exists
$\lambda\in \pbar'(x)\setminus\{1,\tau\}$. For  the same reasons, 
we may assume that $\lambda\not\in \pbar'(F)$ and 
the $\lambda$-sequence  $(w_1, \ldots, w_k)$ is of type B
such that $\pbar'(v_t)=1$. 

 Assume first that the $\lambda$-sequence  is a subsequence of the $\tau$-sequence.  
 Let $w_1=v_i$ for some $i\in [2,t]$, i.e., 
 $(w_1, \ldots, w_k) = (v_i, \ldots, v_t)$ and $\lambda=\pbar'(v_{i-1})$. Since 
 $r\in P_{s_1}(\lambda, \Delta, \phiv')=P_x(\lambda,\Delta, \phiv')$, 
 we do a $(\lambda, \Delta)$-swap at $v_{i-1}$ and get a coloring $\phiv''$.  Since $\{\lambda, \D\}\cap \{1, \tau\} =\emptyset$, 
 $\phiv''$ is also $\{1,\tau\}$-avoiding. 
 But,  under coloring $\phiv''$, the $\tau$-sequence  $(v_1, \ldots, v_{i-1})$ is of type B and
 $\pbar''(v_{i-1})=\Delta$, giving a 
 contradiction to Lemma~\ref{lem:tau-seq2} (v). Under coloring $\phiv'$, for the same reason, the $\tau$-sequence is not a subsequence of the $\lambda$-sequence. 
 
 Assume these two sequences do not have inclusion relation but have a common vertex.  Let $v_i =w_j$ be the first
 common vertex. Clearly, $i\in [2, t]$ and $j\in [2, k]$, and $\pbar'(v_{i-1}) = \pbar'(w_{j-1})=\gamma$. Assume, without loss of generality, $v_{i-1}\notin P_{s_1}(\tau, \D)$, which containing $r$.  We do $(\gamma, \D)$-swap at $v_{i-1}$ and get an $F$-stable coloring $\phiv''$. Since $\{ \gamma, \D\}\cap \{1, \tau\} =\emptyset$,  we have $1, \tau \in \pbar''(x)$ and $\phiv''(s_1u) = \phiv'(s_1u) =1$. But, the $\tau$ sequence $(v_1, \dots, v_{i-1})$ is type B and $\pbar''(v_{i-1}) = \D \ne 1$, a contradiction.

Assume finally  that these two sequences are disjoint. Then one of $v_t$ and $w_k$ are $(1,2)$-unlinked with $x$. Assume, without loss of generality, that  $v_t$ and $r$ are $(1,2)$-unlinked with $x$. 
We do a $(1,2)$-swap at $v_t$ and get a coloring $\phiv''$. Since $s_1u\in P_{s_1}(1, 2) =P_r(1, 2)$, $\phiv''$ is $F$-stable and we still have  
$1, \tau\in \pbar''(x)$ and $\phiv''(s_1u) =1$. But, the $\tau$-sequence $(v_1, \dots, v_t)$ is of type B with 
$\pbar''(v_t) =2 \ne 1$, giving a contradiction. \qed

\section{Proof of Theorem~\ref{thm:longk2}}

We introduce some new concepts in order to prove Theorem~\ref{thm:longk2}. 
\subsection{Pseudo-fan}

  

Let $G$ be a class 2 graph and $rs_1$ be a critical edge. 
A  \emph{pseudo-fan} (P-fan) at $r$ w.r.t. $rs_1$ and a coloring $\varphi\in \CC^k(G-rs_1)$ is a sequence $$S=S_\varphi(r,s_1:s_t:s_p)=(r, rs_1, s_1, rs_2, s_2, \ldots,rs_t, s_t, rs_{t+1},  s_{t+1}, \ldots, s_{p-1},  rs_p, s_p)$$ 
such that  all $s_1, \ldots , s_p$ are distinct vertices in $N_{\D-1}(r)$ and  the following conditions hold:
\begin{enumerate}[(P1)]
	\item $(r, rs_1, s_1, rs_2, s_2, \ldots,rs_t, s_t)$, denoted by $F_\varphi(r,s_1:s_t)$, is a maximum multifan at $r$.
	\item The vertex set $V(S)$ is $\varphi'$-elementary for every $F$-stable $\varphi'$ w.r.t.  $\varphi$.
\end{enumerate}
Clearly every maximum multifan is a P-fan, and if $S$ is a P-fan
w.r.t. $\varphi$ and  $F=F_{\phiv}(r, s_1:s_t)$,  then by the definition above,  $S$ is also a P-fan w.r.t. every   $F$-stable coloring $\varphi'$. The result below is a modification of Lemma~3.6 from~\cite{HZ}. 

\begin{LEM}\label{pseudo-fan-ele}
	Let $G$ be a  class 2 graph with maximum degree $\Delta$, 
	$r\in V_{\Delta}$ be light, and $S=S_\varphi(r, s_1: s_p: s_{q})$ be a P-fan w.r.t. $rs_1$ and a coloring
	$\phiv\in \CC(G-rs_1)$. Then the following two statements hold, where 
	$F=F_{\phiv}(r, s_1:s_p)$.  
	
\begin{enumerate}[(a)]
\item For every $v_1\in V(S)\setminus V(F)$, the $\phiv(rv_1)$-sequence $(v_1,\ldots, v_t)$ is a rotation at $r$,   and  $v_i$   and $r$ are  $(1, \pbar(v_i))$-linked for each $i\in[1,t]$. \label{pseudo-a}

\item For every $i, j$ with $i\in [1,p]$ and $j\in [p+1, q]$ and colors $\gamma\in \pbar(s_i)$ and $\delta\in \pbar(s_j)$, 
$r\in P_{s_i}(\gamma, \delta) = P_{s_j}(\gamma, \delta)$. Moreover, if $\phiv(rz)=\gamma$ for some $z\in N(r)$, then $P_{s_i}(\gamma, \delta)$ meets $z$ before $r$. \label{pseudo-b}

		\label{pseudo-c}
	\end{enumerate}
\end{LEM}
\pf By relabeling colors and vertices, we assume $F$ is typical. Let  $F=F_\varphi(r,s_1:s_\alpha:s_\beta)$ be a typical multifan, where  $\beta=p$.

For   statement \eqref{pseudo-a},  we let $v_1\in V(S)\setminus V(F)$, and  let $\varphi(rv_1)=\tau$.
Since $F$ is maximum, by Lemma~\ref{lem:tau-seq}, we let $(v_1,\ldots, v_t)$ be the $\tau$-sequence at $r$. 
If the sequence is of type A, then it is a rotation and we are done. 
So we assume the sequence is of type B or C. 

If  $(v_1,\ldots, v_t)$ is of type B, i.e., $\pbar(v_t)=\gamma\in \pbar(F)$, since $\pbar_F^{-1}(\gamma)$
and $r$ are $(1,\gamma)$-linked,  
we do a $(1,\gamma)$-swap at $v_t$ to get $\varphi'$. Then we do the B-shifting from $v_1$
to $v_t$ and exchanging the role of 1 and $\tau$ in the entire graph.
This results in an $F$-stable coloring, yet $V(S)$ is not elementary, contradicting (P2)
of the definition of a P-fan.

If  $(v_1,\ldots, v_t)$ is of type C, i.e., $\pbar(v_t)=\pbar(v_{i-1})=\tau_i$ for some $i\in [2,t-1]$ and some $\tau_i\in [1,\D]\setminus \pbar(F)$, since one of $v_{i-1}$ and $v_t$ is 
 $(1,\tau_i)$-unlinked with $r$, 
 we do a $(1,\tau_i)$-swap at a vertex in $\{v_{i-1}, v_t\}$ that is  $(1,\tau_i)$-unlinked with $r$.
 This gives an $F$-stable coloring such that the corresponding $\tau$-sequence is of type B,
 converting the problem to the previous case. 
 
 Thus the $\tau$-sequence  $(v_1,\ldots, v_t)$ is a rotation. Moreover,  $v_i$ and $r$ are $(1,\pbar(v_i))$-linked for each $i\in [1,t]$. As otherwise, a $(1,\pbar(v_i))$-swap at $v_i$ would give rise to a 
 $B$-type $\tau$-sequence, contradicting what  was proved above.  
The proof of statement \eqref{pseudo-a} is completed.

By Statement~\eqref{pseudo-a}, we let $(v_1,\ldots, v_t)$ be the rotation containing $s_j$, where $v_1=s_j$. 
For the first part of statement\eqref{pseudo-b}, suppose to the contrary that
 $r\in P_{s_i}(\gamma, \delta)=P_{v_1}(\gamma, \delta)$ does not hold. Assume without loss of generality that $i\in [1,\alpha]$. Then we have the following three cases: $r\notin P_{s_i}(\gamma, \delta)$ and $r\notin P_{v_1}(\gamma, \delta)$; $r\notin P_{s_i}(\gamma, \delta)$ and $r\in P_{v_1}(\gamma, \delta)$; and $r\in P_{s_i}(\gamma, \delta)$ and $r\notin P_{v_1}(\gamma, \delta)$.

Suppose that $r\notin P_{s_i}(\gamma, \delta)$ and $r\notin P_{v_1}(\gamma, \delta)$. Then let $\varphi'=\varphi/Q$,  where $Q$ is the $(\gamma, \delta)$-chain containing $r$. Note that $\varphi'$ and $\varphi$ agree on every edge incident to $r$ except two edges $rv_2$ and $rz$ where $z$ is the vertex in $N(r)$
such that $\varphi(rz)=\gamma$. Since $r\notin P_{s_i}(\gamma, \delta)$, $r\notin P_{v_1}(\gamma, \delta)$ and $V(S)$ is $\varphi$-elementary, $\pbar'(s_i)=\pbar(s_i)$ for all $s_i\in V(S)$. Thus under the new coloring $\varphi'$, $F^*=(r, rs_1, s_1, \ldots,s_i, rv_2, v_2, \ldots, rv_t, v_t, rv_1, v_1, rs_{i+1},  s_{i+1}, \ldots,s_\beta)$ is a multifan. This is because, if $i<\alpha$,  then $\pbar'(s_i)=\gamma=\varphi'(rv_2)$ and   $\pbar'(v_1)=\delta=\varphi'(rs_{i+1})$, and if $i=\alpha$, then  $\varphi'(s_{i+1})=\Delta \in \pbar'(s_{1})$.  As $|V(F)|<| V(F^*)|$, 
we obtain a contradiction to the maximality of $F$.

Suppose that $r\notin P_{s_i}(\gamma, \delta)$ and $r\in P_{v_1}(\gamma, \delta)$. Then let $\varphi'=\varphi/P_{v_1}(\gamma, \delta)$. Similar to the case above, one can easily check that $F^*=(r, rs_1, s_1, \ldots,s_i, rv_2, v_2, \ldots, rv_t, v_t, rv_1, v_1)$ is a multifan. 
Since $\pbar'(s_i)=\pbar'(v_1)=\gamma$, we obtain a contradiction to 
Lemma~\ref{thm:vizing-fan1}  that 
$V(F^*)$ is $\varphi'$-elementary. 

Suppose that $r\in P_{s_i}(\gamma, \delta)$ and $r\notin P_{v_1}(\gamma, \delta)$. Then let $\varphi'=\varphi/P_{v_1}(\gamma, \delta)$. Note that $\varphi'$ is $F$-stable w.r.t. $\varphi$, thus by the definition of a P-fan, $V(S)$ is $\varphi'$-elementary. But $\pbar'(s_i)=\pbar'(v_1)=\gamma$, giving a contradiction. This completes the proof of the first part of statement \eqref{pseudo-b}.

For the second  part of statement\eqref{pseudo-b}, assume to the contrary that $P_{s_i}(\gamma, \delta)$
meets $r$ before  $z$. Then $P_{s_i}(\gamma, \delta)$ meets $v_2$ before $r$. Let $\varphi'$ be obtained from $\varphi$ by shifting from $v_1$
to $v_t$. Then  $r\not\in P_{s_i}(\delta,\gamma, \varphi')$, showing a contradiction to 
the first part of~\eqref{pseudo-b}. 
\qed 

\subsection{Two structural lemmas}

\begin{LEM}\label{lem:pseodo-adj}
	Let $G$ be a class 2 graph with maximum degree $\Delta \ge 3$, 
	$r\in V_{\Delta}$ be light, and $rs_1$ be a critical edge.   
	If  $S=S(r, s_1:s_p:s_q)$ is a P-fan, 
	then for any $x\in N(S)\setminus N[r]$, $d(x) \ne \D-1$.  
\end{LEM}

\pf Suppose to the contrary that there is a degree $(\D -1)$ vertex $x\notin N[r]$ and a vertex $s^*\in S$ such that  $x\sim s^*$.  
Let $F=F(r, s_1:s_{\alpha}:s_{\beta})$ be the maximum multifan contained in $S$. Since $rs_1$
is a critical edge of $G$,  every edge of $F$ is a critical edge of $G$.  
Thus by Theorem~\ref{thm:s1-adj}, $s^*\in V(S)\setminus V(F)$. 

We may first assume $1\in \pbar(x)$. To see this, 
let $\tau \in \pbar(x)$. If $\tau\in \pbar(F)$, since $\pbar_F^{-1}(\tau)$ and $r$
are $(1,\tau)$-linked by Lemma~\ref{thm:vizing-fan2}~\eqref{thm:vizing-fan1b},
we simply do a $(1,\tau)$-swap at $x$.  Thus we assume that $\tau \in [1,\Delta]\setminus \pbar(F)$. 
By Lemma~\ref{lem:tau-seq2} (i), there is an $F$-stable coloring such that $1$ is missing at $x$.  
We then  do a $(1,\Delta)$-swap at $x$,  
still call it $\varphi$.  We now have  $\Delta \in \pbar(x)$.

We claim that there is an $V(F)$-stable coloring $\phiv'$ such that $\varphi'(s^*x)\in \{2, \Delta\}$. 
Let $\varphi(s^*x)=\tau$. Assume first that $\tau \in \pbar(F)$. If $\tau$ is not $\D$-inducing, we simply do a $(\tau,\Delta)$-swap at $x$. 
Otherwise, we do  $(\Delta,1)-(1,2)-(2,\tau)$-swaps at $x$, and get a desired $V(F)$-stable coloring.  
Thus, we may assume $\tau \in [1,\Delta]\setminus \pbar(F)$. 
For every $V(F)$-stable coloring $\phiv'$, since $F$ is maximum, $r\in P_{s_1}(\tau, \D, \phiv')$ (by Lemma~\ref{lem:rs1}).  We claim
$P_x(\tau, \D, \phiv') = P_{s_1}(\tau, \D, \phiv')$. Otherwise, a $(\tau, \D)$-swap at $x$ gives a desired coloring. 
Applying Lemma~\ref{lem:tau-seq2} (iii), the $\tau$-sequence
$(v_1, \ldots, v_t)$ is of type B and  $\pbar(v_t)=\Delta$. 
Since $r\in P_{s_1}(\tau,\Delta)=P_x(\tau,\Delta)$, we do a $(\tau,\Delta)$-swap at $v_t$ to get an $F$-stable coloring,
and then do the $A$-shifting from $v_1$ to $v_t$. 
Under the new coloring,  $ P_{s_1}(\tau,\Delta) \ne P_x(\tau,\Delta)$. 
Since still $r\in P_{s_1}(\tau,\Delta)$ by Lemma~\ref{lem:rs1}, we do a $(\tau,\Delta)$-swap at $x$
to get a desired $V(F)$-stable coloring.  
So, we  still denote by $\phiv$ and  assume that $\varphi(s^*x)=\Delta$. 

We then show that there is an $V(F)$-stable coloring  $\phiv'$ such that $\varphi'(s^*x)=\Delta$
and $1\in \pbar'(x)$.  
Let $\tau \in \pbar(x)$. If $\tau \in \pbar(V(S))$,  by Lemma~\ref{thm:vizing-fan2}~\eqref{thm:vizing-fan1b} and Lemma~\ref{pseudo-fan-ele}~\eqref{pseudo-a},
we simply do a $(1,\tau)$-swap at $x$. Thus $\tau \in [1,\Delta]\setminus \pbar(S)$. 
If 
there is a  $V(F-r)$-stable  and  $\{ \Delta\}$-avoiding coloring $\varphi'$ such that 
$\pbar'(r)\in \pbar'(x)$, then by switching colors $1$ and $\pbar'(r)$ for the entire graph, we get 
a desired $V(F)$-stable coloring.  So, we assume that there is no such a coloring.  In particular, we have 
$P_x(1, \tau, \phiv) = P_r(1, \tau, \phiv)$.  
By Lemma~\ref{lem:tau-seq2} (ii), 
the $\tau$-sequence  $(v_1, \ldots, v_t)$  at $r$ is of type B such that $\pbar(v_t)=\Delta$. Let $\pbar(s^*)=\delta$. As $V(S)$ is elementary and $\D\in \pbar(s_1)$, we have $v_t\notin V(S)$, and so $s^*\ne v_t$. 
We also note that $\delta \ne \tau$. Otherwise, 
by Lemma~\ref{pseudo-fan-ele}~\eqref{pseudo-a}, $P_{s^*}(1, \tau, \phiv) = P_r(1, \tau, \phiv)$, which gives a contradiction to 
$P_x(1, \tau, \phiv) = P_r(1, \tau, \phiv)$. 
By Lemma~\ref{pseudo-fan-ele}~\eqref{pseudo-b},
$r\in P_{s_1}(\delta,\Delta)=P_{s^*}(\delta,\Delta)$.
We do the $(\delta,\Delta)$-swap at $v_t$ and get a $V(F)$-stable coloring with $\delta$ missing at $v_t$. 
Applying Lemma~\ref{pseudo-fan-ele}~\eqref{pseudo-a} to $s^*\in V(S)$, we get $P_r(1,\delta)=P_{s^*}(1,\delta)$.
We do the $(1,\delta)$-swap at $v_t$. 
Note that  by  Lemma~\ref{pseudo-fan-ele}~\eqref{pseudo-a}, the  $\phiv(rs^*)$-sequence containing $s^*$ at $r$
is a rotation, thus $s^*\not\in \{v_1, \ldots, v_t\}$. 
We do the $B$-shifting from $v_1$ to $v_t$ followed by switching color $1$ and $\tau$ for the entire graph, 
which results in a desired $V(F)$-stable coloring.  

Hence, we may assume that  $\varphi(s^*x)=\Delta$,
$1\in \pbar(x)$, and $\pbar(s^*)=\delta$. 
By Lemma~\ref{pseudo-fan-ele}~\eqref{pseudo-a} that $P_r(1,\delta, \phiv)=P_{s^*}(1,\delta, \phiv)$, 
we do a $(1,\delta)$-swap at $x$.  Under the new coloring, $P_{s^*}(\delta,\Delta)=s_1x$, 
showing a contradiction to that $s^*$ and $s_1$ are $(\delta, \D)$-linked (Lemma~\ref{pseudo-fan-ele}~\eqref{pseudo-b}). 
\qed

\begin{LEM}\label{lem:longk3}
	Let $G$ be a class 2 graph with maximum degree $\Delta \ge 3$, 
	$r\in V_{\Delta-1}$ be light, and $F$ be a multifan at $r$ w.r.t. edge $rs_1$ and 
	a coloring $\phiv   \in \CC^\Delta(G-rs_1)$.
	If $F$ is maximum, then $\pbar(r)\not\subseteq \pbar(x)$ for any $x\in V(G)\setminus N[r]$ with $(N(x)\cap N(s_1))\setminus N_{\Delta-1}[r]\ne \emptyset$.  
	
\end{LEM}

\pf Suppose to the contrary that  there exists a vertex $x\in V(G)\setminus N[r]$ such that $(N(x)\cap N(s_1))\setminus N_{\Delta-1}[r]\ne \emptyset$ and 
$\pbar(r) \subseteq \pbar(x)$.  Let $u \in (N(x)\cap N(s_1))\setminus N_{\Delta-1}[r]$,   $\pbar(r)=\{1,\Delta-1\}$ and $\pbar(s_1)=\{2,\Delta\}$.  So, 
$\{1, \D -1\}\subseteq \pbar(x)$.
Our goal is to modify $\varphi$ in getting a $V(F)$-stable coloring  $\varphi'$ such that $K=(r,rs_1,s_1,s_1u,u,ux,x)$
is a Kierstead path but $\pbar'(x)\cap (\pbar'(s_1)\cup \pbar'(r))\ne \emptyset$, in achieving a contradiction to 
Lemma~\ref{Lemma:kierstead path1}.  Since $r$ is light, we may assume that $F=F(r, s_1:s_{\alpha}:s_{\beta})$ is typical.   

By doing $(1,2)$- and $(\Delta-1,\Delta)$-swaps at $x$ when it is necessary, 
we may assume that $2,\Delta \in \pbar(x)$. 
Applying Lemma~\ref{lem:good-coloring1} (ii)(b), we may assume that 
there is an $V(F)$-stable coloring, still denoted by $\varphi$, such that $\varphi(s_1u)=1$
and $\Delta\in \pbar(x)$. 

We show next that there is an $V(F)$-stable coloring, still denoted by $\varphi$,
such that $\varphi(s_1u)=1$ and $\varphi(ux)=\Delta$. 
Let $\varphi(ux)=\tau$. 

Assume first that $\tau \in \pbar(F)$, if $\tau$
is not $\D$-inducing, we do a $(\tau,\Delta)$-swap at $x$ in getting a desired $V(F)$-stable coloring. 
If $\tau$ is $\Delta$-inducing, we do  $(\Delta,1)-(1,2)-(2,\tau)$-swaps at $x$  in getting a desired 
$V(F)$-stable coloring. 

Suppose that $\tau\in [1,\Delta]\setminus \pbar(F)$.  We claim that for every 
$V(F)$-stable coloring $\phiv'$, $P_x(\tau, \D) = P_{s_1}(\tau, \D)$. Otherwise, 
since $F$ is maximum, $r\in P_{s_1}(\tau,\Delta)$. The $(\tau, \D)$-swap at $x$ gives
a $V(F)$-stable coloring $\phiv'$ such that $\phiv'(ux)=\D$ and $\phiv'(s_1u) =1$, which is what we want. 
By Lemma~\ref{lem:tau-seq2} (vi), the $\tau$-sequence $(v_1, \ldots, v_t)$
is of type B such that $\pbar(v_t)\in \{1,\Delta\}$ or is 2-inducing. 
If $\pbar(v_t)=1$, we do a $(1,2)$-swap at $v_t$, so the color missing at $v_t$ is $2$-inducing. 
Thus we only need to consider two cases. If $\pbar(v_t)=\gamma$ is 2-inducing, 
we do a $(\gamma, \Delta-1)$-swap at $v_t$, where  $\Delta-1$ is another color missing at $r$. 
Then we do the  $B$-shifting from $v_1$ to $v_t$ and get a $V(F-r)$-stable coloring $\phiv'$. In particular, 
we have $\tau\in \pbar'(r)$. Since
$\phiv'(s_1u) =1 \in \pbar'(r)$ and $\phiv'(ux) =\tau \in \pbar'(r)$,  $K=(r,rs_1,s_1,s_1u,u,ux,x)$
is a Kierstead path.  But $\Delta$ is missing at both $s_1$ and $x$, achieving a contradiction to 
Lemma~\ref{Lemma:kierstead path1}. 
Thus, $\pbar(v_t)=\Delta$. Since $r\in P_{s_1}(\tau,\Delta)=P_x(\tau,\Delta)$, we do a $(\tau, \Delta)$-swap at 
$v_t$, resulting in a type A $\tau$-sequence. 
Since  $r\in P_{s_1}(\tau,\Delta)=P_x(\tau,\Delta)$, the $A$-shifting from $v_1$ to $v_t$
gives an $V(F)$-stable coloring such that $x$ and $s_1$ are $(\tau,\Delta)$-unlinked. 
Since still $r\in P_{s_1}(\tau,\Delta)$ by Lemma~\ref{lem:rs1}, we do a $(\tau, \Delta)$-swap at 
$x$ in getting a desired coloring. 

Recall that we assumed that  $\varphi(s_1u)=1$ and $\varphi(ux)=\Delta$. 
Since $\phiv(s_1u) =1 \in \pbar(r)$ and $\phiv(ux) =\D\in \pbar(s_1)$, $K=(r,rs_1,s_1,s_1u,u,ux,x)$
is a Kierstead path. We next show  that there is a $V(F)$-stable coloring $\phiv'$ 
keeping the Kierstead path but 
$\pbar'(x)\cap (\pbar'(s_1)\cup \pbar'(r))\ne \emptyset$, which gives a 
 contradiction to 
Lemma~\ref{Lemma:kierstead path1}.

Let $\tau \in \pbar(x)$. If $\tau \in \pbar(F)$, we simply  do a $(\tau, \Delta-1)$-swap at $x$ 
to get a contradiction. Thus, $\tau \in [1,\Delta]\setminus \pbar(F)$. 
We claim for any $V(F)$-stable coloring $\phiv^*$, $P_x(2, \tau) = P_{s_1}(2, \tau)$. Otherwise, 
since $F$ is  maximum, by Lemma~\ref{lem:rs1} we have $r\in P_{s_1}(2, \tau)$. The $(2, \tau)$-swap at $x$ gives 
a $V(F)$-stable coloring that  maintains the Kierstead path, but $2$ is missing at both $x$ and $s_1$, a contradiction. 
Applying Lemma~\ref{lem:tau-seq2} (vi), the $\tau$-sequence $(v_1, \ldots, v_t)$
is of type B such that $\pbar(v_t)\in \{1,\D\}$ or is 2-inducing. 
If $\pbar(v_t)=1$, we do a $(1,2)$-swap at $v_t$. 
Thus we only need to consider two cases where $\pbar(v_t)\ne 1$. If $\pbar(v_t)=\gamma$ is 2-inducing, 
we do the $(\gamma, \Delta-1)$-swap at $v_t$ and then do the $B$-shifting from $v_1$ to $v_t$. Now $K=(r,rs_1,s_1,s_1u,u,ux,x)$
is a Kierstead path but $\tau$ is missing at both $r$ and $x$, achieving a contradiction to 
Lemma~\ref{Lemma:kierstead path1}. 
Thus, $\pbar(v_t)=\Delta$. Now  doing  $(\Delta,1)-(1,2)-(2,\Delta-1)$-swaps at $v_t$
 and then the $B$-shifting from $v_1$ to $v_t$
gives a same contradiction as right before.  
\qed

\subsection{Proof of Theorem~\ref{thm:longk2}}

Since all vertices not missing  a given color $\alpha$
are saturated by the matching that consists of all edges colored by $\alpha$ in $G$, we have the 
following result.  

\begin{LEM}[Parity Lemma]
	Let $G$ be an $n$-vertex graph and $\varphi\in \CC^\Delta(G)$. 
	Then for any color $\alpha\in [1,\Delta]$, 
	$|\{v\in V(G): \alpha\in \pbar(v)\}| \equiv n \pmod{2}$. 
\end{LEM}

\begin{THM4}\label{thm:longk2a}
	Let $G$ be a  $\Delta$-critical graph with $n$ vertices.  If  $G$ has a light $\D$-vertex  and $\Delta  >n/2+1$,
	then $n$ is odd. 
\end{THM4}

\pf  Let $r$ be a light  $\D$-degree vertex of $G$. Recall that $N(r)=N_\Delta(r)\cup N_{\Delta-1}(r)$. 
We prove first that 
$d(x)=\Delta$ for every $x\in V(G)\setminus N[r]$. 
Assume to the contrary that there exists $x\in V(G)\setminus N[r]$ with $d(x) \le \Delta-1$. If $d(x)\ge \Delta-2 \ge (n-1)/2 $, since as 
 $d(r) = \D \ge (n+3)/2$, we get $|N(r)\cap N(x)| \ge  d(r) + d(x) - |N(x)\cup N(r)| \ge (n+2) - (n-2) =4$. Since $|N_{\D}(r)| =2$,   there exists $s\in N_{\Delta-1}(r)$
such that $x\sim s$. Since $G$ is $\Delta$-critical, $rs$ is a critical edge of $G$. 
But this gives a contradiction to Theorem~\ref{thm:s1-adj}. Thus $d(x)\le \Delta-3$. 
Then for any $u\in N_\Delta(x)$, there  exists $s\in N_{\Delta-1}(r)$
such that $u\sim s$. Since $d(x)\le \Delta-3$ and every neighbor of $r$ has degree at least $\Delta-1$,
we have $u\not\in N(r)\setminus N_{\Delta-1}(r)$. 
Again, using that $rs$ is a critical edge of $G$, we obtain a contradiction to Theorem~\ref{thm:longk}.

Assume to the contrary that $n$ is even. Let $s\in N_{\D-1}(r)$ and $\phiv\in \CC^{\D}(G-rs)$. By the Parity Lemma, if a color is missing at one vertex, it must be missing at another vertex. Let $X$ be a largest $\varphi$-elementary set that contains $r$ and $s$ such that $X\subseteq N_{\Delta-1}[r]$. 
Then since all vertices in $V(G)\setminus N_{\Delta-1}(r)$ are of maximum degree, 
\begin{equation*}
|N_{\Delta-1}(r)\setminus X|=|\pbar(N_{\Delta-1}(r)\setminus X)| \ge |\pbar(X)|=|X|+1. 
\end{equation*}
 On the other hand, we have 
 \begin{equation*}
 |N_{\Delta-1}(r)\setminus X|+|X\setminus\{r\}|=\D-2. 
 \end{equation*}
Combining the two formulas above, we get $|X|\le \frac{\D-2}{2}$. Thus $X$
contains at most $\frac{\D-4}{2}$ vertices from $N_{\D-1}(r)$. 


We claim that for every $s\in N_{\D -1}(r)$, 
$|N(s)\cap N_{\Delta-1}(r)|\le  \frac{\Delta-4}{2} $.  Otherwise,  we let $s$ be such a vertex and let $\varphi\in \CC^\Delta(G-rs)$, and let $X$
be defined the same way as above.  Since $|X| \le \frac{\Delta-4}{2}$ and $|N(s)\cap N_{\D -1}(r)| > \frac{\Delta-4}{2}$, 
there is a vertex $x\in N(s)\cap (N_{\D-1}(r) \setminus X)$. So, the color $\tau =\phiv(rx)$ presents at every vertex. 
Let $G_1$ be obtained from $G$ by deleting all the edges colored by $\tau$.
Then $G_1$ is still a class 2 graph and  $r$ is a light maximum degree vertex in $G_1$, and $\varphi\in \CC^{\D-1}(G_1-rs)$.  
However, in $G_1$, $x\not\sim r$  and $x\sim s$, contradicting Theorem~\ref{thm:s1-adj}. 
Thus we assume that  for every $s\in N_{\Delta-1}(r)$, it holds that 
$|N(s)\cap N_{\Delta-1}(r)|\le \frac{\Delta-4}{2} $. 

Let $N_{\Delta-1}(r)=\{s_1,\ldots, s_{\Delta-2}\}$, 
 $\varphi\in \CC^\Delta(G-rs_1)$, and let $X$ be a largest $\varphi$-elementary set that contains $r$ and $s_1$ such that $X\subseteq N_{\Delta-1}[r]$.
Since $|X|\le \frac{\Delta-4}{2}$ and $|N_{\D-1}(r)| = \D -2$, there exist a vertex
$x\in N_{\D -1}(r)$ such that the color $\tau =\phiv(rx)$  is presented at every vertex of $G$. 
Let $G_1$ be obtained from $G$ by deleting all the edges colored by $\tau$. 
Then $G_1$ is still a class 2 graph such that $r$ is a light maximum degree vertex, and $\varphi\in \CC^{\Delta-1}(G_1-rs_1)$. 
As $\Delta(G_1)=\Delta-1\ge n/2+1$,  there exists $s^*\in N_{G_1}(r)$ with $d_{G_1}(s^*) = \Delta(G_1)-1$
such that $x\sim s^*$ in $G_1$. Note that $G_1$
is still a class 2 graph, and $\varphi$, with being restricted on $G_1$, is a $\Delta(G_1)$-coloring of $G_1$. 
Let $F_\varphi(r,s_1:s_\alpha:s_\beta)$ be a maximum typical multifan at $r$
and $S$ be a maximum P-fan containing $F$. 
If $s^*\in V(S)$, then we obtain a contradiction to Lemma~\ref{lem:pseodo-adj}. 
Thus $s^*\not\in V(S)$. 
Since $V(S)$ is a largest P-fan containing $F$, there is a $V(F)$-stable coloring $\varphi$ such that 
$V(S)\cup \{s^*\}$ is not $\varphi$-elementary. 
Since $V(S)$ is $\varphi$-elementary by the definition of $S$, $\pbar(s^*)\in \pbar(S)$. 
As for every $\gamma \in \pbar(S)\setminus \pbar(r)$,  $\pbar^{-1}_F(\gamma)$ and $r$
are $(1,\gamma)$-linked by Lemma~\ref{thm:vizing-fan2}~\eqref{thm:vizing-fan1b} 
and Lemma~\ref{pseudo-fan-ele}~\eqref{pseudo-a}, we do a $(1,\pbar(s^*))$-swap at $s^*$. 
Let $\varphi(rs^*)=\delta$.

If $rs^*$ is a critical edge of $G_1$, then we already reach a contradiction to Theorem~\ref{thm:s1-adj}. 
Thus, $rs^*$ is not a critical edge of $G_1$. 
We let $G_2=G_1-rs^*$. Note that $G_2$ is still a class 2 graph with $r \in V_{\Delta(G_2)-1}$
being a light vertex.  The coloring  $\varphi$, with being restricted on $G_2$, is a $\Delta(G_2)$-coloring of $G_2$, and  
$F_\varphi(r,s_1:s_\alpha:s_\beta)$ is still  a maximum typical multifan at $r$.
By the choice of $\varphi$ before, we have $\pbar(r)=\pbar(s^*)=\{1,\delta\}$.

Since $s_1$ is adjacent in $G_2$ to at most $\frac{\Delta-4}{2}$ vertices from $\{s_1,\ldots, s_{\Delta-2}\}$, 
and $d_{G_2}(s_1) \ge \Delta-2$,  $s_1$ is adjacent in $G_2$
to at least $\Delta/2-1$ vertices from $V(G)\setminus \{r, s_1,\ldots, s_{\Delta-2}\}$.
Similarly, $d_{G_2}(s^*)=\Delta(G_2)-2= \Delta-3$,  
$s^*$ is adjacent in $G_2$ to at most $\frac{\Delta-4}{2}$ vertices from $\{s_1,\ldots, s_{\Delta-2}\}$,
and $s^*\not\sim r$, 
it follows that $s^*$  is adjacent in $G_2$
to at least $\Delta/2-1$ vertices from $V(G)\setminus \{r, s_1,\ldots, s_{\Delta-2}\}$.
Since $\Delta \ge n/2+2$, $|V(G)\setminus \{r, s_1,\ldots, s_{\Delta-2}\}| \le n/2-1 $.
As  $2(\Delta/2-1) \ge n/2$, there exists $u\in (N_{G_2}(s_1)\cap N_{G_2}(x))\setminus\{r, s_1,\ldots, s_{\Delta-2}\}$. Since $\pbar(r)=\pbar(s^*)=\{1,\delta\}$, we obtain a contradiction to Lemma~\ref{lem:longk3}.  
\qed

\bibliographystyle{plain}
\bibliography{SSL-BIB_08-19}

\begin{thebibliography}{10}

\bibitem{HZ}
Y.~Cao, G.~Chen, G.~Jing, and S.~Shan.
\newblock Proof of the {C}ore {C}onjecture of {H}ilton and {Z}hao.
\newblock {\em \arxiv{2004.00734}}, 2020.

\bibitem{MR848854}
A.~G. Chetwynd and A.~J.~W. Hilton.
\newblock Star multigraphs with three vertices of maximum degree.
\newblock {\em Math. Proc. Cambridge Philos. Soc.}, 100(2):303--317, 1986.

\bibitem{MR975994}
A.~G. Chetwynd and A.~J.~W. Hilton.
\newblock The edge-chromatic class of graphs with maximum degree at least
  {$|V|-3$}.
\newblock In {\em Graph theory in memory of {G}. {A}. {D}irac ({S}andbjerg,
  1985)}, volume~41 of {\em Ann. Discrete Math.}, pages 91--110. North-Holland,
  Amsterdam, 1989.

\bibitem{fw}
S.~Fiorini and R.~J. Wilson.
\newblock Edge-colourings of graphs, {R}esearch notes in {M}aths.
\newblock {\em Pitman, London}, 1977.

\bibitem{MR0349458}
Jean-Claude Fournier.
\newblock Colorations des ar\^{e}tes d'un graphe.
\newblock {\em Cahiers Centre \'{E}tudes Recherche Op\'{e}r.}, 15:311--314,
  1973.
\newblock Colloque sur la Th\'{e}orie des Graphes (Brussels, 1973).

\bibitem{Gupta-67}
R.~G. Gupta.
\newblock {\em Studies in the Theory of Graphs}.
\newblock PhD thesis, Tata Institute of Fundamental Research, Bombay, 1967.

\bibitem{MR1395947}
A.~J.~W. Hilton and Cheng Zhao.
\newblock On the edge-colouring of graphs whose core has maximum degree two.
\newblock {\em J. Combin. Math. Combin. Comput.}, 21:97--108, 1996.

\bibitem{comm}
D.~G. Hoffman and C.~A. Rodger.
\newblock The chromatic index of complete multipartite graphs.
\newblock {\em J. Graph Theory}, 16(2):159--163, 1992.

\bibitem{Holyer}
Ian Holyer.
\newblock The {NP}-completeness of edge-coloring.
\newblock {\em SIAM J. Comput.}, 10(4):718--720, 1981.

\bibitem{K-path-84}
Henry~A. Kierstead.
\newblock On the chromatic index of multigraphs without large triangles.
\newblock {\em J. Combin. Theory Ser. B}, 36(2):156--160, 1984.

\bibitem{MR2082738}
Michael Plantholt.
\newblock Overfull conjecture for graphs with high minimum degree.
\newblock {\em J. Graph Theory}, 47(2):73--80, 2004.

\bibitem{seymour79}
P.~D. Seymour.
\newblock On multicolourings of cubic graphs, and conjectures of {F}ulkerson
  and {T}utte.
\newblock {\em Proc. London Math. Soc. (3)}, 38(3):423--460, 1979.

\bibitem{StiebSTF-Book}
Michael Stiebitz, Diego Scheide, Bjarne Toft, and Lene~M. Favrholdt.
\newblock {\em Graph edge coloring}.
\newblock Wiley Series in Discrete Mathematics and Optimization. John Wiley \&
  Sons, Inc., Hoboken, NJ, 2012.
\newblock Vizing's theorem and Goldberg's conjecture, With a preface by
  Stiebitz and Toft.

\bibitem{Vizing64}
V.~G. Vizing.
\newblock On an estimate of the chromatic class of a {$p$}-graph.
\newblock {\em Diskret. Analiz}, (3):25--30, 1964.

\bibitem{vizing-2factor}
V.~G. Vizing.
\newblock The chromatic class of a multigraph.
\newblock {\em Kibernetika (Kiev)}, 1965(3):29--39, 1965.

\bibitem{Vizing-2-classes}
V.~G. Vizing.
\newblock Critical graphs with given chromatic class.
\newblock {\em Diskret. Analiz No.}, 5:9--17, 1965.

\end{thebibliography}
\end{document}